\def\hal{1}

\title{Event-triggered output feedback stabilization via dynamic high-gain scaling}
\documentclass[journal]{IEEEtran}
\IEEEoverridecommandlockouts  
\usepackage{amsfonts,amsmath,mathtools}
\usepackage{graphicx,color}
\usepackage{url}
\usepackage{hyperref} 

\input{newcommands_mathGeneral}
\newcommand{\alphaMax}{\alpha ^\ast}
\newcommand{\deltaMax}{\delta ^\ast}
\newcommand{\deltaMaxo}{\delta_o ^\ast}
\newcommand{\deltaMaxc}{\delta_c ^\ast}

\newcommand{\LmatC}{\mathcal{L}}

\newcommand{\Sc}{\mathcal{S}}

\newcommand{\Next}{_{k+1}}
\newcommand{\NextM}{_{k+1}^-}
\renewcommand{\L}{L_k}
\newcommand{\LM}{L_k^-}
\newcommand{\Lnext}{L_{k+1}}
\newcommand{\LnextM}{L_{k+1}^-}	

\newcommand{\Lmat}{\mathcal{L}_k}

\newcommand{\xNext}{x_{k+1}}

\newcommand{\xo}{\hat x_k}
\newcommand{\xNexto}{\hat x_{k+1}}
\newcommand{\xNextMo}{\hat x_{k+1}^-}

\newcommand{\Xo}{\hat X_k}
\newcommand{\Xnexto}{\hat X_{k+1}}
\newcommand{\XnextMo}{\hat X_{k+1}^-}

\newcommand{\SnextM}{\mathcal{S}_{k+1}^-}
\newcommand{\Snext}{\mathcal{S}_{k+1}}
\renewcommand{\S}{\mathcal{S}_{k}}
\providecommand{\e}{}
\renewcommand{\e}{e_k}
\newcommand{\eNext}{e_{k+1}}
\newcommand{\eNextM}{e_{k+1}^-}

\newcommand{\E}{E_k}
\newcommand{\Enext}{E_{k+1}}
\newcommand{\EnextM}{E_{k+1}^-}

\def\figurename{Fig.}
\def\tablename{Tab.}

\newcommand{\MatlabSource}[1]{\href{#1}{\underline{sources Matlab}}}

\definecolor{souris}{gray}{0.65}

\newtheorem{theorem}{Theorem}

\newtheorem{lemma}{Lemma}
\newtheorem{proposition}{Proposition}

\newtheorem{assumption}{Assumption}

\newtheorem{remark}{Remark}

\newcommand{\Figure}[3]{
  \begin{figure}[!ht]
    \centering
    \includegraphics[width=#2\linewidth]{#1}
    \caption{#3}\label{fig:#1}
  \end{figure}
}
\newcommand{\FIG}[1]{\figurename~\ref{#1}}

\renewcommand{\|}{|}
\newcommand{\Fc}{F_c}
\newcommand{\Fo}{F_o}
\newcommand{\Ro}{R_o}
\newcommand{\Rc}{R_c}
\newcommand{\dr}{\mathfrak d}

\catcode`\@=11
\def\downparenfill{$\m@th\braceld\leaders\vrule\hfill\bracerd$}
\def\overparen#1{\mathop{\vbox{\ialign{##\crcr\crcr
\noalign{\kern0.4ex}
\downparenfill\crcr\noalign{\kern0.4ex\nointerlineskip}
$\hfil\displaystyle{#1}\hfil$\crcr}}}\limits}
\catcode`\@=12

\author{Johan Peralez, Vincent Andrieu, Madiha Nadri, Ulysse Serres
\thanks{All authors are with Universit\'e Lyon 1 CNRS UMR 5007 LAGEP, France.
(e-mail johan.peralez@gmail.com, vincent.andrieu@gmail.com, nadri@lagep-lyon1.fr, ulysse.serres@univ-lyon1.fr)
}
\thanks{V.~Andrieu is also with  Fachbereich C - Mathematik und Naturwissenschaften, Bergische Universit\"at Wuppertal, Germany.
}
\thanks{This work was supported by ANR LIMICOS contract number 12 BS03 005 01.}}

\begin{document} 

\maketitle

\begin{abstract}
This work addresses output feedback stabilization via event triggered output feedback.
In the first part of the paper, linear systems are considered,
whereas the  second part shows that a dynamic event triggered output feedback control law can achieve feedback stabilization of the origin for a class of nonlinear systems by employing dynamic high-gain techniques.
\end{abstract}

\section{Introduction}
The implementation of a control law on a process requires the use of an appropriate sampling scheme.
In this regards, periodic control (with a constant sampling period) is the usual approach that is followed for practical implementation on digital platforms.
Indeed, periodic control benefits from a huge literature, providing a mature theoretical background (see e.g. \cite{Astroem1997,Nesic1999,alur2007}) 
and numerous practical examples.
The use of a constant sampling period makes  closed-loop analysis and implementation easier,
allowing solid theoretical results and a wide deployment in the industry.
However, the rate of control execution being fixed by a worst case analysis (the chosen period must guarantee the stability for all possible operating conditions), this may lead to an unnecessary fast sampling rate and then to an overconsumption of available resources. 

The recent growth of shared networked control systems for
which communication and energy resources are often limited goes with an increasing interest in aperiodic control design.
This can be observed in the comprehensive overview on event-triggered and self-triggered control presented in \cite{Heemels2014}.
Event-triggered control strategies introduce a triggering condition assuming a continuous monitoring of the plant (that requires a dedicated hardware) 
while in self-triggered strategies, the control update time is based on predictions using previously received data. The main drawback of self-triggered control is the difficulty to guarantee an acceptable degree of robustness, especially in the case of uncertain systems.  

Most of the existing results on event-triggered and self-triggered control for nonlinear systems are based on the input-to-state stability (ISS) assumption which implies the existence of a feedback control law ensuring an ISS property with respect to measurement errors 
(\cite{tabuada2007event,Anta2010,Abdelrahim2015a,Postoyan2015}) and also \cite{seuret2013stability}.

In this ISS framework, an emulation approach is followed: the knowledge of an existing robust feedback law in continuous time is assumed,
and some triggering conditions are proposed to preserve stability under sampling.

Another proposed approach consists in the redesign of a continuous time stabilizing control. 
For instance,
the authors in \cite{Marchand2013} adapted the original \textit{universal formula} introduced by Sontag for nonlinear control affine systems. 
The relevance of this method was experimentally shown in \cite{Villarreal-Cervantes2015} where the regulation of an omnidirectional mobile robot was addressed. 

Although aperiodic control literature has demonstrated an interesting potential,
important fields still need to be further investigated to allow a wider practical deployment.
In particular,
literature on output feedback control for nonlinear systems  is scarce (\cite{Yu2013}, \cite{Abdelrahim2014}, \cite{Liu2015}, \cite{Tanwani2015CDC}) whereas, in many control applications, the full state information is not available for measurement.

The high-gain approach is a very efficient tool to address the stabilizing control problem in the continuous time case. It has the advantage to allow uncertainties in the model and to remain simple. 

Different approaches based on high-gain techniques have been followed in the literature to tackle the output feedback problem in the continuous-time case (see for instance \cite{Andrieu2009,Krishnamurthy2004,andrieu2007unifying,AndrieuTarbouriech_TAC_2013}) and more recently for the (periodic) discrete-in-time case (see \cite{Qian2012}). 
In the context of observer design, \cite{Andrieu2014} proposed the design of a continuous discrete time observer, revisiting high-gain techniques in order to give an adaptive sampling stepsize (see also \cite{dinh2015continuous,mazenc2015design} for observers with constant sampling period). 

In this work,
we extend the results obtained in \cite{Andrieu2014} to event-triggered output feedback control. 
In high-gain designs, the asymptotic convergence is obtained by dominating the nonlinearities with high-gain techniques. 
In the proposed approach,  high-gain is dynamically adapted with respect to time varying nonlinearities in order to allow an efficient trade-off between the high-gain parameter and the sampling step size. Moreover, the proposed strategy is shown to ensure the existence of a minimum inter-execution time.
Note that a preliminary version of this work has appeared in \cite{Peralez2016Nolcos} in which only an event triggered state feedback was considered.

The paper is organized as follows.
The control problem and the class of considered systems is given in Section \ref{sec_PrbStatement}.
In Section \ref{sec_PrelimResult}, some preliminary results concerning linear system are given.
The main result is stated in Section \ref{Sec_MainResult} and its proof is given in Section \ref{sec:stateFeedback_analysis}.
Finally Section \ref{sec_IllustrativeExample} contains an illustrative example.

\section{Problem Statement} \label{sec_PrbStatement}
\subsection{Class of considered systems}
In this work, we consider the problem of designing an event-triggered output feedback for the class of uncertain nonlinear systems described by the dynamical system  
\begin{align} 
	\dot{x}(t) &= Ax(t) + Bu(t) + f(x(t)), \label{eq:problem_statement_sys}
\end{align}
where the state $x$ is in $\RR^n$; $u:\RR\rightarrow\RR$ is the control signal in
$\LL^\infty(\RR_+,\RR)$, $A$ is a matrix in $\RR^{n\times n}$ and $B$ is a vector in $\RR^{n}$ in the following form
\begin{equation} \label{eq:problem_statement_A_B_f}
A=\begin{bmatrix} 0 & 1 & 0 & \cdots & 0 \\ \vdots & \ddots & \ddots & \ddots & \vdots \\ 0 & \cdots & 0 & 1 & 0 \\ 0 & \cdots & \cdots &0 & 1 \\ 0 & \cdots & \cdots &\cdots & 0\end{bmatrix} , \qquad 
B=\begin{bmatrix} 0 \\ \vdots \\ 0 \\ 0  \\ 1 \end{bmatrix},\end{equation}
 and $f:\RR^n\rightarrow\RR^n$ is a vector field having the following triangular structure
\begin{equation}
f(x) = \begin{bmatrix} f_1(x_1) \\ f_2(x_1,x_2) \\ \vdots \\ f_n(x_1,x_2, \dots, x_n) \end{bmatrix}.
\end{equation}

We consider the case in which the vector field $f$ satisfies the following assumption.
\begin{assumption}[Nonlinear bound] \label{hyp:incremental_bound}
There exist a non-negative continuous function $c$,  
positive real numbers $c_0$, $c_1$ and $q$ such that 
for all $x\in \RR^n$, we have 
\begin{align}\label{E:assumption1}
	|f_j(x(t))| \leq& c(x_1) \left(|x_1| + |x_2| + \dots + |x_j| \right),	
\end{align}	 
with 
\begin{align}	
	c(x_1) =& c_0 + c_1 |x_1|^q.
\end{align}	 
\end{assumption}
Notice that Assumption~\ref{hyp:incremental_bound} is more general than the incremental property
introduced in \cite{Qian2012}  since the function $c$ is not constant but depends on $x_1$.
This bound can be also related to \cite{praly2003asymptotic,Krishnamurthy2004} in which continuous output feedback laws were designed.
Note however that in these works no bounds were imposed on the function $c$. Moreover, in our present context we do not consider inverse dynamics.

\subsection{Updated sampling time controller}

In the sequel, we restrict ourselves to a  sample-and-hold implementation, i.e.
the input is assumed to be constant between any two execution times.
The control input $u$ is defined through a sequence $(t_k,u_k)_{k\in\NN}$ in $\RR_+\times\RR$
in the following way 
\begin{equation}\label{eq_SampledControl}
	u(t) = u_k, \quad \forall~t \in [t_k, t_{k+1})\ .
\end{equation}
It can be noticed that for $u$ to be well defined for all positive time, we need that
\begin{equation}\label{eq_WellDefinedSampling}
\lim_{k\rightarrow+\infty} t_k = +\infty\ .
\end{equation}

Our control objective is to design the sequence   $(u_k,t_k)_{k\in\NN}$ such that the origin of the obtained closed loop system is asymptotically stable.
%
This sequence depends only on the output which in our considered model is simply given as
\begin{equation}\label{eq_Output}
y(t) = Cx(t)\ ,\ C = \begin{bmatrix}
1 & 0 & \cdots & 0
\end{bmatrix}\ .
\end{equation}

\Figure{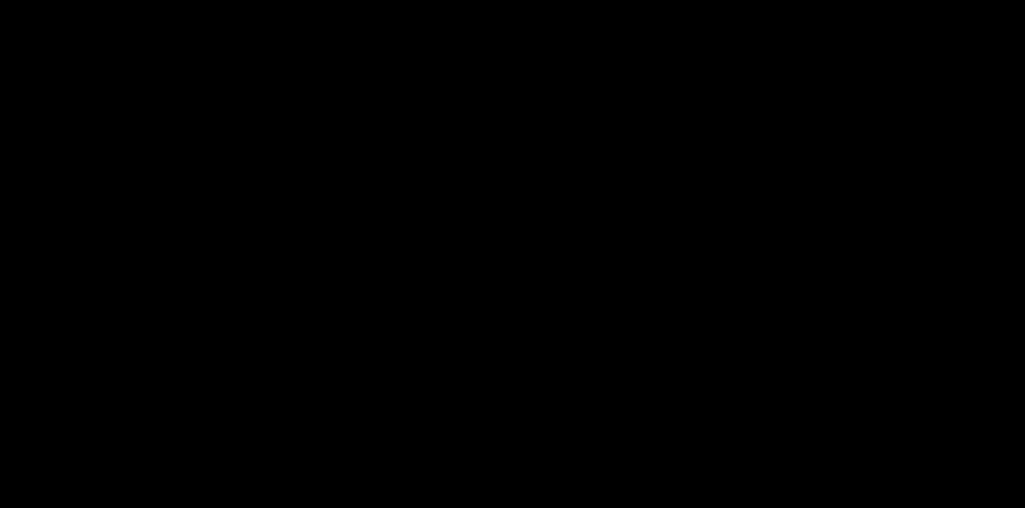}{1}{Event-triggered control schematic.}

Note however that in the same spirit as for  the sample and hold control, we consider only
a sequence of output values
\begin{equation}\label{eq_Output_Sampled}
y_k = Cx(t_k)\ ,
\end{equation}
which corresponds to the evaluation of the output  $y(\cdot)$ at the same time instant $t_k$.

In addition to a feedback controller that
computes the control input, event-triggered and self-triggered control systems
need a \textit{triggering mechanism} that
determines when a new measurement occurs and when the control input has to be updated again. This rule is said to be \textit{static} if it only involves the current state of the system, and \textit{dynamic} if it uses an additional internal dynamic variable \cite{Girard2015}.
Our approach is summarized  in  \FIG{fig:schema_Event}.
	

\subsection{Notation}
In this paper, we denote by $\langle\cdot,\cdot\rangle$ the canonical scalar product in $\RR^n$ and by $\|\cdot\|$ the induced Euclidean norm; we use the same notation for the corresponding induced matrix norm. 
Also, we use the symbol $'$ to denote the transposition operation.

To simplify the presentation, we introduce the following notations:
$\xi(t^-) = \lim\limits_{\substack{\tau \to t \\ \tau <t}}\xi(\tau)$,  $\xi_k = \xi(t_k)$ and $\xi_k^- = \xi(t_k^-)$.

\section{Preliminary result: linear case}\label{sec_PrelimResult}
In high-gain design, the idea is to consider the nonlinear terms (the $f_i$'s) as disturbances. 
A first step consists in synthesizing a robust control for the linear part of the system, neglecting the effects of the nonlinearities. Then, convergence and robustness are amplified through a high gain parameter to deal with the nonlinearities.

Therefore, let us first focus on a general linear dynamical system
\begin{align} \label{eq:sys_linear_case}
	\dot{x}(t) &= A x(t) + B u(t),
\end{align}
where the state $x$ evolves in $\RR^n$ and the control $u$ is in $\RR$. The matrix $A$ is in $\RR^{n\times n}$ and the matrix $B$ is in $\RR^{n}$.
The measured output is given as a sequence of values $(y_k)_{k\geq0}$ in $\RR$ as in \eqref{eq_Output_Sampled}
where $C$ is a column vector in $\RR^{n}$ and $(t_k)_{k\geq0}$ is a sequence of times to be selected.

In this preliminary case, we review a well known result concerning periodic sampling approaches. Indeed, an emulation approach is adopted for the stabilization of the linear part: a feedback law is designed in continuous time and a triggering
condition is chosen to preserve stability under sampling. 

It is well known that if there exists a continuous time dynamical output feedback control law that asymptotically stabilizes the system,
then there exists a positive  inter-execution time $\delta=t\Next-t_k$ such that the sampled control law renders the system asymptotically stable.  This result is rephrased in the following Lemma \ref{lemm:trigg_linear} whose proof 
is postponed to Appendix \ref{sec:proof_lemma_linear}.
\begin{lemma} \label{lemm:trigg_linear}
\renewcommand{\O}{(A+BK)}
Suppose that there exist a row vector $K_c$ and a column vector $K_o$ (both in $\RR^{n}$) rendering $(A+BK_c)$ and $(A+K_o C)$ Hurwitz.
Then there exists a positive real number $\deltaMax$ such that for all $\delta$ in $[0; \deltaMax)$ the following holds.
Let the sequence
$(t_k,u_k)_{k\in\NN}$ be defined as
\begin{equation}\label{eq:th_trigg_linear_Lk_1}
 t_0 = 0\ ,\ t_{k+1}=t_k+ \delta \ ,\ 
 u_k =  K_c \hat x(t_k)\ ,\ \forall \ k\in\NN\ ,
\end{equation}
where $\hat x(t_0)$ is in $\RR^n$ and for $k$ in $\NN^*$ 
\begin{align}
\label{eq:lin_SystObs2}
\dot {\hat x}(t)
&= A\hat x(t) + Bu_k , \quad\! \forall~  t\in \left[t_k,t_{k+1}\right), \\
\hat x(t_k)
&= \hat x(t_k^-) + \delta   K_o (C\hat x(t_k^-)- y_k). \label{eq:lin_SystObs}
\end{align}
Then $(x(t),\hat x(t))=0$ is a globally and asymptotically stable (GAS) solution  for  the dynamical system
defined by
\eqref{eq_SampledControl}, \eqref{eq:sys_linear_case}, \eqref{eq:th_trigg_linear_Lk_1}, \eqref{eq:lin_SystObs2} and \eqref{eq:lin_SystObs}.
\end{lemma}

\medskip
This result which is based on robustness is valid for general matrices $A$, $B$ and $C$.

We want to point out that the proof of Lemma \ref{lemm:trigg_linear} is based on the fact that  if $A+BK_c$ and $A+K_oC$ are Hurwitz, the origin of the discrete time linear system defined for all $k$ in $\NN$ as
\begin{equation}\label{eq_Discrete}
\begin{bmatrix}
\xNexto\\
\eNext
\end{bmatrix} = 
\begin{bmatrix}\displaystyle
 \Fc(\delta)  &  \delta K_o C \exp(A \delta)\\
 0 &  \Fo(\delta)
\end{bmatrix}\begin{bmatrix}
\xo\\
e_k
\end{bmatrix}
\end{equation}
where $e=\hat x - x$  is the estimation error,  and
\begin{align}
\Fc(\delta) &= \exp(A \delta) + \int_{0}^{\delta} \exp(A(\delta-s))B K_c ds \label{eq_Fc}\\
\Fo(\delta) &= (I + \delta K_oC) \exp(A \delta)\label{eq_Fo}
\end{align}
is asymptotically stable for $\delta$ sufficiently small.

However, when we consider the particular case in which $(A,B,C)$ are as in (\ref{eq:problem_statement_A_B_f}) and (\ref{eq_Output}) 
(i.e. an integrator chain),  it is shown in the following theorem that the inter-execution time can be selected arbitrarily large as long as the control is modified.
\begin{theorem}[Chain of integrator] \label{th:trigg_linear}
Suppose the matrices $A$, $B$ and $C$ have the structure stated in \eqref{eq:problem_statement_A_B_f}-\eqref{eq_Output}.
Let  $K_c$ and $K_o$ both in $\RR^n$, be such that $A+BK_c$ and $A+K_oC$ are Hurwitz.
Then there exists a positive real number $\alphaMax$ such that for all $\alpha$ in $[0,\alphaMax)$ the following holds.\\
For all $\delta>0$, let the sequence
$(t_k,u_k)_{k\in\NN}$ be defined as
\begin{equation}\label{eq:th_trigg_linear_Lk}
 t_0 = 0\ ,\ t_{k+1}=t_k+ \delta \ ,\ 
 u_k =  K_cL^{n+1}\mathcal{L} \hat x(t_k)\ ,\ 
\end{equation}
where  $\hat x(t_0)$ is in $\RR^n$ and for $k$ in $\NN^*$ 
\begin{align}
\label{eq:lin_SystObs3}
\dot {\hat x}(t)
&= A\hat x(t) + Bu_k , \quad\! \forall~  t\in \left[t_k,t_{k+1}\right), \\
\hat x(t_k)
&= \hat x(t_k^-) + \delta \mathcal{L}^{-1}  K_o (C\hat x(t_k^-)- y_k), \label{eq:lin_SystObs2_jump}
\end{align}
and
\begin{equation}\label{eq_alphaLdelta}
\mathcal{L} =\diag\left(\frac{1}{L}, \dots,\frac{1}{L^n}\right)\ ,\ L = \dfrac{\alpha}{\delta}\ .
\end{equation}
Then $(x(t),\hat x(t))=0$ is a GAS solution  for  
the dynamical system
defined by
\eqref{eq_SampledControl}, \eqref{eq:sys_linear_case}, \eqref{eq:th_trigg_linear_Lk}, \eqref{eq:lin_SystObs3} and \eqref{eq:lin_SystObs2_jump}.
\end{theorem}

\begin{remark}
 Note that the difference between equation \eqref{eq:lin_SystObs} and equation \eqref{eq:lin_SystObs2_jump} is the $\mathcal{L}^{-1}$ factor that appears in the latter.
\end{remark}

\begin{remark}
 Note that in the particular case of the chain of integrator the sampling period $\delta$ can be selected arbitrarily large.
 To obtain this result the two \textit{gains} $K_c$ and $K_o$ have to be modified as seen in equations \eqref{eq:th_trigg_linear_Lk} and \eqref{eq:lin_SystObs2_jump}
\end{remark}

\begin{proof}
In order to analyze the behavior of the closed-loop system, let us mention the following algebraic properties of the matrix $\mathcal{L}$:
\begin{equation}\label{eq:linear_algebraic_prop}
\mathcal{L} A = L A \mathcal{L}\ , \  \mathcal{L}B  = \frac{B }{L^n}\ ,\ C\LmatC^{-1}=LC\ .
\end{equation}
Let $e=\hat x - x$.
Consider now the  following change of coordinates
\begin{equation}\label{eq:ScaledCoordinate}
\hat X	= \mathcal L \hat x 
\ ,\ E	= \mathcal L  e 
\end{equation}
Employing \eqref{eq:linear_algebraic_prop} and \eqref{eq:th_trigg_linear_Lk}, it yields that in the new coordinates the closed-loop dynamics are
for all $t$ in $[t_k, t_{k+1})$:
\begin{align}
\dot {\hat X}(t) =& L \left( A\hat X(t) + BK_c\Xo \right),		\\
\dot E(t) =& L A E(t).
\label{eq:trigg_linear_sys_new_coord}
\end{align}
By integrating the previous equality and  employing $L\delta = \alpha$, it yields for all $k$ in $\NN$:
\begin{align*}
\XnextMo &=  \left[\exp(AL\delta) + \int_0^\delta\exp(AL(\delta-s))LBK_cds\right]\Xo \\
		 &=  \Fc(\alpha)\Xo	,\\
\EnextM 		&= \exp(A\alpha)\E,		 	
\end{align*}
and with \eqref{eq:lin_SystObs2_jump}
\begin{align*}
\Xnexto 
		&=\mathcal L \left( \xNextMo + \delta \mathcal L^{-1} K_o  C \eNextM \right)\\
		&= \XnextMo + \alpha K_o C \EnextM	\\
		&= \Fc(\alpha)\Xo 	+ \alpha K_o C \exp(A\alpha)\E\ .
\end{align*}
Similarly, it yields:
\begin{align*}			
\Enext 
		&= \mathcal L (I + \delta \mathcal L^{-1} K_o  C) \eNextM\\
		&= (I+ \alpha K_o C) \EnextM \\
		&= \Fo(\alpha)\E\ .
\end{align*}
In other words,
this is the same discrete dynamic as the one given in \eqref{eq_Discrete}.
Consequently,
from Lemma \ref{lemm:trigg_linear}, there exists a positive real number $\alphaMax$ such that $(\hat X,E)=0$ (and thus $(x,\hat x)=0$) is a GAS equilibrium for the system \eqref{eq:trigg_linear_sys_new_coord} provided $L \delta$ is in $[0, \alphaMax)$. 
\end{proof}

\section{Main result: the nonlinear case}
\label{Sec_MainResult}
We now consider the full nonlinear system \eqref{eq:problem_statement_sys} with $f$ satisfying Assumption \ref{hyp:incremental_bound}.
Following the high-gain paradigm, the considered control law is the one used for the chain of integrator in \eqref{eq:th_trigg_linear_Lk}-\eqref{eq:lin_SystObs3}-\eqref{eq:lin_SystObs2_jump} with \eqref{eq_SampledControl}. 
In the context of a linear growth condition, i.e. if the bound $c(x_1)$ defined in Assumption \ref{hyp:incremental_bound} is replaced by a constant $c$,
the authors have shown in \cite{Qian2012} that a (well chosen) constant parameter $L$ can guarantee the global stability, provided that $L$ is greater than a function of the bound.
However,
with a bound in the form \eqref{E:assumption1} of Assumption \ref{hyp:incremental_bound}, we need to adapt the high-gain parameter to follow a function of the time varying bound.
Following the idea presented in \cite{Andrieu2014} in the context of observer design, we define $L$  as the evaluation at time $t_k^-$ of the following continuous discrete dynamics:
\begin{align}
	\dot{L}(t) &= a_2 L(t)M(t)c(x_1(t)), & \forall t \in [t_k,t_k + \delta_k) 	 \label{eq:high_gain_update_law_start}	\\	
	\dot{M}(t) &= a_3 M(t) c(x_1(t)), & \forall t \in [t_k,t_k + \delta_k) \label{eq:high_gain_update_law_M} \\
	\L &= \LM (1-a_1 \alpha)+ a_1 \alpha \label{eq:high_gain_update_law_Lnext} \\
	M_{k} &= 1, \label{eq:high_gain_update_law_end}
\end{align}
with initial condition $L(0) \geq 1$, $M(0)= 1$ and where $a_1, a_2 ,a_3$ are positive real numbers to be chosen. For a justification of this type of high-gain update law, the interested reader may refer to \cite{Andrieu2014} where it is shown that this update law is a continuous discrete version of the high-gain parameter update law introduced in \cite{praly2003asymptotic}.

With this high-gain parameter and following what has been done in Theorem \ref{th:trigg_linear}, the sequence of control is defined as follows.
\begin{equation} \label{eq:control_law_start}
u_k=   K_cL_k^{n+1}\mathcal{L}_k \hat x(t_k)\ , \ \forall k\in \NN,
\end{equation}
where $\hat x(0)$ is in $\RR^n$.
And, for $k$ in $\NN^*$
\begin{align}
\dot {\hat x}(t)	=& A\hat x(t) + Bu(t) , \quad\! \forall~  t\in \left[t_k,t_{k+1}\right), \label{eq:control_law_obs}\\
\hat x(t_k)	=& \hat x(t_k^-) + \delta_{k-1} (\mathcal{L}^-_k)^{-1}   K_o (C\hat x(t_k^-)- y_k).\label{eq:control_law_obs_jump}
\end{align}
with
$\mathcal{L}_k^- =\diag\left(\frac{1}{L_k^-}, \dots,\frac{1}{(L_k^-)^n}\right)$.

It remains to select the sequences  $\delta_k$  and the execution times $t_k$ .
These are given by the following relations,
\begin{align}
	t_0 = 0, \quad  t_{k+1}=t_k + \delta_k, \label{eq:control_law_tNext}
\end{align}
\begin{equation}
\delta_k = \min\{s\in \RR_+\ |\  s L((t_k+s)^-)=\alpha\}. \label{eq:control_law_end}
\end{equation}
Equations \eqref{eq:control_law_tNext}-\eqref{eq:control_law_end} constitute the triggering mechanism of the self-triggered strategy. 
It does not directly involve the state value $x$ but the additional dynamic variable $L$ and so can be referred as a dynamic triggering mechanism (\cite{Girard2015}). The relationship between $\L$ and $\delta_k$ comes from the right hand side equation of \eqref{eq_alphaLdelta}. It highlights the trade-off between high-gain value and inter-execution time (see \cite{Dabroom2001,Qian2012}).

We are now ready to state our main result whose proof is given in Section \ref{sec:stateFeedback_analysis}.

\begin{theorem}
\emph{(Stabization via event-triggered output feedback control):} 
\label{th:main_result}
Assume the  functions $f_i$'s in \eqref{eq:problem_statement_sys}   satisfy Assumption~\ref{hyp:incremental_bound}.
Then, there exist positive numbers $a_1$, $a_2$, $a_3$, two gain matrices $K_c$, $K_o$ and $\alphaMax>0$ such that for all $\alpha$ in $[0,\alphaMax]$,
there exists a positive real number $L_{\max}$  such that the set
\[
\{x=0,\hat x=0, L\leq L_{\max}\}\subset \RR^n\times \RR^n \times \RR,
\]
is GAS along the solution of system \eqref{eq:problem_statement_sys} with  
the self-triggered feedback  \eqref{eq:control_law_start}-\eqref{eq:control_law_end}.
More precisely, there exists a class $\mathcal{KL}$ function $\beta$ such that by denoting $(x(\cdot),\hat x(\cdot), L(\cdot))$ the solution initiated from $(x(0),\hat x(0), L(0))$ with  $L(0)\geq 1$, this solution is defined for all $t\geq 0$ and satisfies
\begin{multline}\label{eq_MainResult}
|x(t)| + |\hat x(t)| + |\tilde L(t)| 
\\\leq 
\beta(|x(0)| + |\hat x(0)| + |\tilde L(0)| ,t),
\end{multline}
where $\tilde L(t) = \max\{L(t)-L_{\max}\}$.
Moreover there exists a positive real number $\delta_{\min}$ such that $\delta_k>\delta_{\min}$ for all $k$ and so ensures the existence of
a minimal inter-execution time.
\end{theorem}


\section{Proof of Theorem \ref{th:main_result}} \label{sec:stateFeedback_analysis}
Following \cite{praly2003asymptotic}, let us introduce the following scaled coordinates along a trajectory of system \eqref{eq:problem_statement_sys}
which will be used at different places in this paper (compare with \eqref{eq:ScaledCoordinate}).
\begin{equation}\label{eq:ScaledCoordinate_b}
\hat X(t)	= \mathcal S(t) \hat x(t)
\ ,\ E(t)	= \mathcal S(t)  e(t)\ , 
\end{equation}
where
$$
\mathcal S(t)  =L(t)^{1-b}\mathcal L(t)\ ,\ \mathcal L(t)=\diag\left(\frac{1}{L(t)}, \dots,\frac{1}{L^n(t)}\right),
$$
 $e(t)=\hat x(t) - x(t)$, and 
where $1\geq b> 0$ is such that $bq<1$ with $q$ given in Assumption~\ref{hyp:incremental_bound}.
%

\renewcommand{\O}{(A+BK)}
\subsection{Selection of the gain matrices $K_c$ and $K_o$}
Let $D$ be the diagonal matrix in $\RR^{n\times n}$ defined by $D=\diag(b,1+b,\dots, n+b-1)$. 
Let $P$ and $Q$ be two symmetric positive definite matrices and $K_c$, $K_o$ two vectors in $\RR^n$ such that (always possible, see \cite{andrieu2009high}) 
\begin{align} 
		P(A+BK_c) + (A+BK_c)'P &\leq -I, \label{eq:bound_K_Hurwitz}\\
		p_1 I\leq P &\leq p_2 I, \label{eq:bound_P} \\
		p_3 P \leq PD+DP_c &\leq p_4 P,\label{eq:bound_D}\\
		Q(A+K_oC) + (A+K_oC)'Q & \leq -I, \label{eq:bound_Ko_Hurwitz}\\
		q_{1} I \leq Q & \leq q_{2} I, \label{eq:bound_Po} \\
		q_{3} Q \leq QD+DQ & \leq q_{4} Q, \label{eq:bound_Do}
\end{align}
with $p_1, \dots, p_4, q_{1}, \dots, q_{4}$ positive real numbers. 

With the matrices $K_c$ and $K_o$ selected it remains to select the parameters $a_1$, $a_2$, $a_3$ and $\alpha^*$.\\
This is done on two steps: in Proposition \ref{pr:existenceXk} we focus on the existence of the sequence $(x_k,L_k)$ for all $k$ in $\NN$. Then,  Proposition \ref{pr:decreaseVk} shows using a Lyapunov analysis that a sequence of quadratic function of scaled coordinates is decreasing.

Based on these two propositions, the proof of Theorem  \ref{th:main_result} is given in Section \ref{Sec_ProofTh} where it is  shown that the time function $L$ satisfies an ISS property (see Proposition \ref{Prop_BoundL}).

\subsection{Existence of the sequence $(t_k,\hat x_k,e_k,L_k)_{k\in\NN}$}

The first step of the proof is to show that the sequence $(\hat x_k,e_k, L_k)_{k\in\NN} = (\hat x(t_k),e(t_k), L(t_k))_{k\in\NN}$ is well defined. 
Note that it does not imply that $(\hat x(t),e(t))$ is defined for all $t$ since for the time being it has not been shown that the sequence $t_k$ is unbounded.
This will be obtained in Section \ref{Sec_ProofTh} when proving Theorem \ref{th:main_result}.

\begin{proposition}[Existence of the sequence] \label{pr:existenceXk}
Let $a_1$, $a_3$ and $\alpha$ be positive, and $a_2\geq \frac{3n}{q_1}$, where $q_1$ was defined in (\ref{eq:bound_Po}).
Then, the sequence $(t_k,\xo, \e, \L)_{k\in\NN}$ is well defined.
\end{proposition}
\noindent\textbf{Proof of Proposition \ref{pr:existenceXk}:}
We proceed by contradiction.
Assume that $k\in\NN$ is such that
$(t_k,\xo, \e, \L)$ is well defined but $(t_{k+1},\xNexto, \eNext, \Lnext)$ is not.
This means that there exists a time $t^*>t_k$
such that $\hat x(\cdot)$, $e(\cdot)$ and $L(\cdot)$ are well defined for all $t$ in $[t_k,t^*)$ and such that 
\begin{equation}\label{eq_assunbd}
	\lim_{t\rightarrow t^*}\big( |\hat x(t)| + |e(t)| +|L(t)| \big) = +\infty.
\end{equation}
Since
$L(\cdot)$ is increasing and, in addition, for all $t$ in $[t_k, t^*)$ we have (according to (\ref{eq:control_law_end})) $L(t) \leq \frac{\alpha}{(t-t_k)}$,
we get:
\begin{equation}\label{eq:BoundL}
	L^*=\lim_{t\rightarrow t^*}L(t)\leq \frac{\alpha}{(t^*-t_k)}<+\infty.
\end{equation}
Consequently, $\lim_{t\rightarrow t^*} |\hat x(t)| + |e(t)| = +\infty$,
which together with \eqref{eq:ScaledCoordinate_b} yields
\begin{equation}\label{eq_assunbd_X}
	\lim_{t\rightarrow t^*} |\hat X(t)| + |E(t)| = +\infty.
\end{equation}
On the other hand, let $U$ and $W$ be the two quadratic functions
\begin{equation}\label{eq_defUW}
U(\hat X) = \hat X'P\hat X\ ,\ W(E) = E'QE .
\end{equation} 
With a slight abuse of notation, when evaluating these functions along the solution of \eqref{eq:problem_statement_sys}, we denote
$U(t) = U(\hat X(t))$ and $W(t) = W(E(t))$. 
For all $t$ in $[t_k, t^*)$, we have
\begin{align}
	\dot{U}(t)
		&= \dot{\hat X}(t)'P \hat X(t) + \hat X(t)'P\dot{\hat X}(t), \label{eq:dot_V} \\
\dot W(t)
		&= \dot{E}(t)'Q E(t) + E(t)'Q\dot{E}(t),   \label{eq:dot_Vo}
\end{align}
where
\begin{align*}
	\dot{\hat X}(t) 
    &= \dot{\Sc}(t)\hat x(t) + \Sc(t) \dot{\hat x}(t) ,\\
    &= - \dfrac{\dot{L}(t)}{L(t)}D \hat X(t) + L(t) A \hat X(t) + L(t) B K  \Xo ,
\end{align*}
and
\begin{align*}
	\dot{E}(t) 
    &= \dot{\Sc}(t)E(t) + \Sc(t) \dot{E}(t) ,\\
    &= - \dfrac{\dot{L}(t)}{L(t)}D E(t) + L(t) A E(t) - \Sc(t) f\left(x(t)\right).
\end{align*}
With the previous equalities, \eqref{eq:dot_V}-\eqref{eq:dot_Vo} become for all $t$ in $[t_k, t^*)$
\begin{align*}
\dot U(t) &= - \dfrac{\dot L(t)}{L(t)} \hat X(t)'(PD +DP) \hat X(t) \\ & + L(t) [\hat X(t)' (A'P + PA) \hat X(t) + 2 \hat X(t)'P BK \Xo],   \\
	 \dot W(t) &= - \dfrac{\dot L(t)}{L(t)} E(t)'(QD +DQ) E(t) \\ & + L(t) E(t)' (A'Q + QA) E(t) + 2E(t)'Q\Sc(t)f(x(t)).
\end{align*}
Since $M \ge 1$, we have with \eqref{eq:high_gain_update_law_start}, \eqref{eq:bound_D} and \eqref{eq:bound_Do} for all $t$ in $[t_k, t^*)$
\begin{align*}
-\dfrac{\dot L(t)}{L(t)} \hat X(t)'(PD +DP)\hat X(t) 
	\leq& -p_3 a_2 c(x_1(t)) U(t), \\
-\dfrac{\dot L(t)}{L(t)} E(t)'(QD +DQ)E(t) 
	\leq& -q_3 a_2 c(x_1(t))   W(t).    
\end{align*}
Moreover, using Young's inequality, we get
\begin{align*}
	2\hat X(t)'P BK \Xo 
    	&\leq \hat X(t)'P \hat X(t) + \Xo' (K'B'P+PBK) \Xo .
\end{align*}
Hence, taking $\lambda_1$ and $\lambda_2$ such that
$$
A'P + PA + I \leq \lambda_1 P\ ,\ K'B'P+PBK\leq \lambda_2 P\ ,
$$ we have, for all $t$ in $[t_k, t^*)$
\begin{equation}\label{eq_Udot}
\dot U(t)\leq \left(- p_3 a_2 c(x_1(t))  + L(t) \lambda_1 \right)  U(t)+ L(t)\lambda_2 U_k.
\end{equation}
On another hand, with Assumption \ref{hyp:incremental_bound} and since $L(t)\geq 1$,
it yields
\begin{align}
\|\Sc(t)f(x(t))\|^2 &= \sum_{i=1}^n \left|\frac{f_i(x(t))}{L(t)^{i+b-1}}\right|^2 ,\nonumber\\
&\leq  \sum_{i=1}^n \left(c(x_1(t))\sum_{j=1}^i |X_j(t)|\right)^2 ,\nonumber\\
&\leq n^2 c(x_1(t))^2 \|\hat X(t)-E(t)\|^2 .\label{eq_SF}
\end{align}
Hence, we get
\begin{multline*}
2E(t)'Q\Sc(t)f(x(t)) \\\leq 2nc(x_1(t))q_3 \left(\frac{3}{2}E(t)' E(t) + \hat X(t)' \hat X(t)\right).
\end{multline*}
Taking $\lambda_3$  such that
$
A'Q + QA \leq \lambda_3 Q$ and since $2n q_3 I\leq \frac{2nq_3}{p_1}P$
it yields
\begin{multline}\label{eq_Wdot}
\dot W(t) \leq \left(\left(\frac{3n}{q_1}- a_2\right)q_3 c(x_1(t)) + L(t) \lambda_3 \right)W(t) \\
+\frac{2nq_3}{p_1} c(x_1(t))U(t) .
\end{multline}
Let us denote
\begin{equation}\label{eq_LyapV}
V(t) =  U(t) + \mu W(t)\ ,
\end{equation}
where $\mu$ is any positive real number that will be useful in the proof of Proposition \ref{pr:decreaseVk}.
Bearing in mind that $L(t) \leq L^*$ for all $t$  in $[t_k, t^*)$ (from (\ref{eq:BoundL})) and
with the couple $(a_2,\mu)$ selected to satisfy
$a_2 \geq \frac{3n}{q_{1}}  $ and   $a_2p_3\geq  \mu\lambda_4 $, 
 inequalities \eqref{eq_Udot} and \eqref{eq_Wdot} yield 
\begin{align*}
	\dot V(t)
	&\leq L^* \lambda_1 U(t)+ L^* \lambda_2 U_k +  \mu  L^* \lambda_3  W(t), \\ 
    &\leq L^* (\lambda_1 + \lambda_3) V(t) + L^* \lambda_2 V_k .    
\end{align*}
This with \eqref{eq:BoundL} give for all $t$ in $[t_k,t^*)$
\begin{align}
	V(t)
	&\le 
    \begin{multlined}[t]
    	 \exp\left((\lambda_1+\lambda_3) L^* (t-t_k) \right)V_k \\
	     +\int_0^{t-t_k} \exp\big( (\lambda_1+\lambda_3)  L^*  (t-t_k-s)\big) \lambda_2 V_k ds
     \end{multlined}\notag \\
	&\le k(\alpha) V_k,
   	\label{eq:InterSamp}
\end{align}
where $k(\alpha) = \exp\left((\lambda_1+\lambda_3) \alpha \right) +
        (\exp( (\lambda_1+\lambda_3)\alpha)-1) \frac{\lambda_2 }{\lambda_1+\lambda_3}$.
Hence, $\lim_{t\rightarrow t^*} |E(t)| + |\hat X(t)|<+\infty$
which contradicts (\ref{eq_assunbd_X}) and thus, ends the proof.
\null\hfill$\Box$

\subsection{Lyapunov analysis} \label{sec:Lyapunov_analysis}
The second step of the proof of Theorem  \ref{th:main_result} consists in a Lyapunov analysis to show that a good
selection of the parameters $a_1$, $a_2$ and $a_3$ in the high-gain update law \eqref{eq:high_gain_update_law_start}-\eqref{eq:high_gain_update_law_end} yields the decrease of the sequences $V_k=V(t_k)$ defined from \eqref{eq_LyapV} with a proper selection of $\mu$. 


\begin{remark}
Using the results obtained in \cite{praly2003asymptotic} on lower triangular systems, the dynamic scaling \eqref{eq:ScaledCoordinate_b} includes a number $b$. Although the decreases of $V_k$ can be obtained with $b=1$, it will be required that $bq<1$ in order to ensure the boundedness of $L(\cdot)$ (see equation (\ref{eq_usebq}) in Section \ref{sec:L_boundedness}).
\end{remark}

The aim of this subsection is to show the following intermediate result.
\begin{proposition}[Decrease of scaled coordinates] \label{pr:decreaseVk}
There exist $a_1>0$ (sufficiently small), $a_2>0$ (sufficiently large), a continuous function $N$ and $\alphaMax>0$ such that
for $a_3=2n$ and for all $\alpha$ in $[0,\alphaMax]$ there exists $\mu$  such that with the time function $V$ defined in \eqref{eq_LyapV} the following property is satisfied:
\begin{equation}\label{eq:LyapDiscrete}
	V_{k+1} -V_k
	\leq -\alpha N(\alpha) \left(\frac{L_k}{L_{k+1}^-}\right)^{2(n-1+b)} V_k.
\end{equation}
\end{proposition}
\vspace{0.5em}\noindent\textbf{Proof of Proposition \ref{pr:decreaseVk}:}
First of all, we assume that $a_2\geq \frac{3n}{q_1}$. Hence, with Proposition \ref{pr:existenceXk}, we know that the sequence $(t_k,x_k, e_k,L_k)$ is well defined for all $k$ in $\NN$.
Let $k$ be in $\NN$.
The nonlinear system \eqref{eq:problem_statement_sys} with the control \eqref{eq:control_law_start} gives the closed-loop dynamics
\begin{align*}
\dot {\hat x}(t)
&= A\hat x(t) + BK_c(L_k)^{n+1}\mathcal{L}_k \hat x(t_k) ,\\ 
\dot{e}(t)
&= Ae(t)- f(x(t)), &\forall~t \in [t_k,t_k + \delta_k).
\end{align*}
Integrating the preceding equalities between $t_k$ and $t\NextM$ yields
\begin{align}
	\xNextMo
	&=
    \begin{multlined}[t]
    	 \exp(A \delta_k)\xo \\
	     + \int_{0}^{\delta_k} \exp(A(\delta_k-s))B K_cL_k^{n+1}\mathcal{L}_k \xo ds,
     \end{multlined}\notag \\
	\eNextM
	&=
	\exp(A \delta_k)e_k - \int_{0}^{\delta_k} \exp(A (\delta_k-s)) f (x(s))ds, 
   	\label{eq:outputFeed_eNextM}
\end{align}
and with \eqref{eq:control_law_obs_jump}, we get
\begin{align}
	\xNexto
	&=
    \begin{multlined}[t]
    	\exp(A \delta_k)\xo + \delta_k (\mathcal{L}^-_{k+1})^{-1} K_o C \eNextM \\
	    + \int_{0}^{\delta_k} \exp(A(\delta_k-s))B K_cL_k^{n+1}\mathcal{L}_k \xo ds \\
     \end{multlined}\notag \\
	\eNext
	&=
    \begin{multlined}[t]
    	(I + \delta_k (\mathcal{L}^-_{k+1})^{-1} K_oC) \Bigg(\exp(A \delta_k)e_k \\
	    - \int_{0}^{\delta_k} \exp(A (\delta_k-s)) f (x(s))ds \Bigg).  \label{eq:outputFeed_eNext}
     \end{multlined}
\end{align}
In the following, we successively consider the evolution of the $e$ part of the dynamics and the evolution of $\hat x$ part.

\medskip
\noindent\emph{\textbf{Analysis of the term in $e$:} }
Employing the algebraic equality given in \eqref{eq:linear_algebraic_prop} yields that $\LmatC \exp(As) =\exp(L As) \LmatC$.
Hence,
when left multiplying \eqref{eq:outputFeed_eNext} by $\mathcal S_{k+1}^-=\left(L_{k+1}^-\right)^{1-b}\mathcal L_{k+1}^-$, we get the following inequality:
$$
\mathcal S_{k+1}^- \eNext = \Fo (\alpha)\mathcal S_{k+1}^- e_k + \Ro,
$$
where we have used the notations
\begin{align*}
\Fo (\alpha) &= (I + \alpha  K_oC)\exp(A \alpha) ,\\
\Ro  &= - (I + \alpha  K_oC)\int_{0}^{\delta_k} \exp(L_{k+1}^-A (\delta_k-s)) \mathcal S_{k+1}^-f (x(s))ds .
\end{align*}
Let $W_k=W(\E)$ where $W$ and $\E$ are respectively defined in \eqref{eq_defUW}  and  \eqref{eq:ScaledCoordinate_b}. 
Note that, since we have $\Enext = \Psi \SnextM \eNext$ with $\Psi=\Snext \inv{\SnextM}$, it yields from \eqref{eq:outputFeed_eNext}
$$
W_{k+1}
	= W( \Psi \SnextM \eNext) 
	= W_k +T_{o,1} + T_{o,2} ,
$$
with
\begin{align*}
T_{o,1} =& W( \Psi \Fo (\alpha)\mathcal S_{k+1}^- e_k)  - W(\E) ,\\
T_{o,2} =& 2 e_k' \SnextM \Fo (\alpha)' \Psi P \Psi \Ro  + \Ro '\Psi P \Psi \Ro .
\end{align*}
Let $\beta$ be defined by
\[
\beta = n\int_{0}^{\delta_k} c(x_1(t_k+s))ds.
\]
The following two lemmas are devoted to upper bound the two terms $T_{o,1}$ and $T_{o,2}$.
The term $T_{o,1}$ will be shown to be negative thanks to \cite[Lemma 3, p109]{Andrieu2014} which in our context becomes the following Lemma.
\begin{lemma}[\cite{Andrieu2014}] \label{lemm:outputFeed_boundTo1}
Let $a_1\leq \frac{1}{2q_{2}q_{4}}$ and $a_3=2n$.
There exists $\alpha_{o}^*>0$ sufficiently small  such that for all $\alpha$ in $[0,\alpha_{o}^*)$
\begin{align}	
T_{o,1} \le
	-\Bigg(
	\frac{a_2q_{3} q_{1} }{a_3} \left[e^{ 2\beta }-1\right] 
 + \frac{\alpha q_{1} }{4q_{2} }\Bigg) 
\norm{ \SnextM e_k } ^2. \label{eq:outputFeed_boundTo1}
\end{align}
\end{lemma}
For the second term, we have the following estimate.
\begin{lemma}\label{lemm:state_Feed_To2} 
There exist two positive real valued continuous functions $N_{o,\hat x}$ and $N_{o,e}$ such that the following inequality holds
\begin{multline*}
T_{o,2} 	
		\leq 
\left[e^{ 2\beta }-1\right] \left[N_{o,\hat x}(\alpha)  \norm{\SnextM \xo}^2
\right.\\ + \left. 
N_{o,e}(\alpha)  \norm{\SnextM \e}^2\right].
\end{multline*}
\end{lemma}
The proof of Lemma \ref{lemm:state_Feed_To2}
is postponed to Appendix \ref{sec:Proof_Lemm_To2}.

\noindent\emph{\textbf{Analysis of the term in $\hat x$}:}
Employing the algebraic equality given in \eqref{eq:linear_algebraic_prop}, we get
from \eqref{eq:outputFeed_eNext} 
$$
 \xNext = \inv{\mathcal S_{k}}\Fc (\alpha_k) \hat X_k + \inv{\SnextM}\Rc ,
$$
where $\Fc$ is defined in \eqref{eq_Fc}, $\alpha_k = \delta_kL_k$ and
$$
\Rc =\alpha  K_o C E_{k+1}^-.
$$
Let $U_k=U(\hat X_k)$ where $U$ is defined in \eqref{eq_defUW}.
This yields with the former equality
\begin{align*}
U_{k+1}
	=& \xNext'\mathcal S_{k+1}P\mathcal S_{k+1} \xNext \\	
	=& U_k +T_{c,1} + T_{c,2},
\end{align*}
with
\begin{align}
T_{c,1} =&  U(\mathcal S_{k+1} \inv{\mathcal S_{k} }\Fc (\alpha_k) \hat X_k)-U_k,\\
T_{c,2}	=& 2 \Rc' \Psi P\mathcal S_{k+1} \inv{\mathcal S_{k} }\Fc (\alpha_k) \hat X_k + U(\Psi \Rc). \notag 
\end{align}
Similarly, the following two lemmas are devoted to upper bound the two terms $T_{c,1}$ and $T_{c,2}$.
The first one is \cite[Lemma 5.4]{Peralez2015} which is devoted to upper bound $T_{c,1}$.
\begin{lemma}[\cite{Peralez2015}]
\label{lemm:outputFeed_boundTc1}
Let $a_1\leq  \frac{2}{p_4p_2}$ and $a_3=2n$.
Then, there exists $\alphaMax>0$ sufficiently small  such that for all $\alpha$ in $[0,\alphaMax)$
\begin{equation} \label{eq:state_feed_boundT1}
T_{c,1}
\leq -\left(\frac{\alpha}{p_2}\right)^2 U(\hat X_k)
-\| \SnextM \hat x_k \|^2
 (e^{2\beta}-1)\dfrac{p_3 p_1 a_2}{2n} .
\end{equation}
\end{lemma}
The proof of Lemma \ref{lemm:outputFeed_boundTc1} 
can be found in \cite{Peralez2015}.
\begin{lemma}\label{lemm:state_Feed_Tc2}
There exist three positive real valued continuous functions $N_{c,\hat x}$, $N_{c,e}$ and $N_{c,0}$ such that the following inequality holds
\begin{multline*}
T_{c,2} \leq N_{c,1}(\alpha)\norm{\SnextM e_k}^2  + \frac{1}{2}\left(\frac{\alpha}{p_2}\right)^2 U(\hat X_k)\\ 
+\left[N_{c,\hat x}(\alpha)\norm{\SnextM \xo}^2 + N_{c,e}(\alpha)\norm{\SnextM e_k}^2\right] [e^{2\beta}-1] 
\end{multline*}
\end{lemma}
The proof of Lemma \ref{lemm:state_Feed_Tc2} 
is postponed to Appendix \ref{sec:Proof_Lemm_T2}.

\medskip
\noindent\emph{\textbf{End of the proof of Poposition \ref{pr:decreaseVk} }:}
Let $\alpha^*=\max\left\{\alpha_{o}^*,\alpha_{c}^*\right\}$ and let $0<\alpha<\alpha^*$, $a_1 = \min\left\{\frac{1}{2q_{2}q_{4}}, \frac{2}{p_4p_2}\right\}$ and $a_3=2n$.
With Lemma \ref{lemm:outputFeed_boundTo1}, Lemma \ref{lemm:state_Feed_To2}, Lemma \ref{lemm:outputFeed_boundTc1}
 and Lemma \ref{lemm:state_Feed_Tc2},  it yields
\begin{multline}
V_{k+1} - V_k \leq 
-\frac{1}{2}\left(\frac{\alpha}{p_2}\right)^2U_k
+\left[N_{c,1}(\alpha) -	 \mu\frac{\alpha q_{1} }{4q_{2} }\right]\norm{ \SnextM e_k } ^2
\\+\left[e^{ 2\beta }-1\right]  \left[\mu N_{o,\hat x}(\alpha)-\dfrac{p_3 p_1 a_2}{2n} \right]  \norm{\SnextM \xo}^2
\\+
\left[e^{ 2\beta }-1\right]\left[\mu N_{o,e}(\alpha)-\mu \frac{a_2q_{3} q_{1} }{a_3}\right] \norm{\SnextM \e}^2 .
\end{multline}
Taking $\mu$ sufficiently large such that
$$
N_{c,1}(\alpha) -	 \mu\frac{\alpha q_{1} }{4q_{2} }\leq -\frac{1}{2}\mu\frac{\alpha q_{1} }{q_{2} },
$$
 and then $a_2$ sufficiently large such that,
 $$
 \mu N_{o,\hat x}(\alpha)-\dfrac{p_3 p_1 a_2}{2n}\leq 0\ ,\ \mu N_{o,e}(\alpha)-\mu \frac{a_2q_{3} q_{1} }{a_3}\leq 0,
 $$
 it yields
$$
	V_{k+1} -V_k
	\leq -\alpha N_0(\alpha)\left[ U_k + \norm{ \SnextM e_k } ^2 \right] .
$$
where $N_0$ is a continuous function taking postiive values.
Employing the fact that $\frac{L_k}{ L_{k+1}^-}\leq 1$, it yields
$$
\norm{ \SnextM e_k } ^2=\norm{ \SnextM \inv{\S}E_k } ^2 \geq \left(\frac{L_k}{L_{k+1}^-}\right)^{2(n-1+b)}\frac{U_k}{p_2},
$$
which gives the existence of a continuous function $N$ such that inequality \eqref{eq:LyapDiscrete} holds. 
This ends the proof of Proposition \ref{pr:decreaseVk}.

\begin{remark}
Due to the jumps of the high-gain parameter $L$ at instants $t_k$ in equation \eqref{eq:high_gain_update_law_Lnext}, the Lyapunov function $t\mapsto V(t)$ does not decrease continuously as illustrated in \FIG{fig:V_principe}.
However,
the sequence $(V_k)_{k\geq0}$ is decreasing.
\end{remark}

\Figure{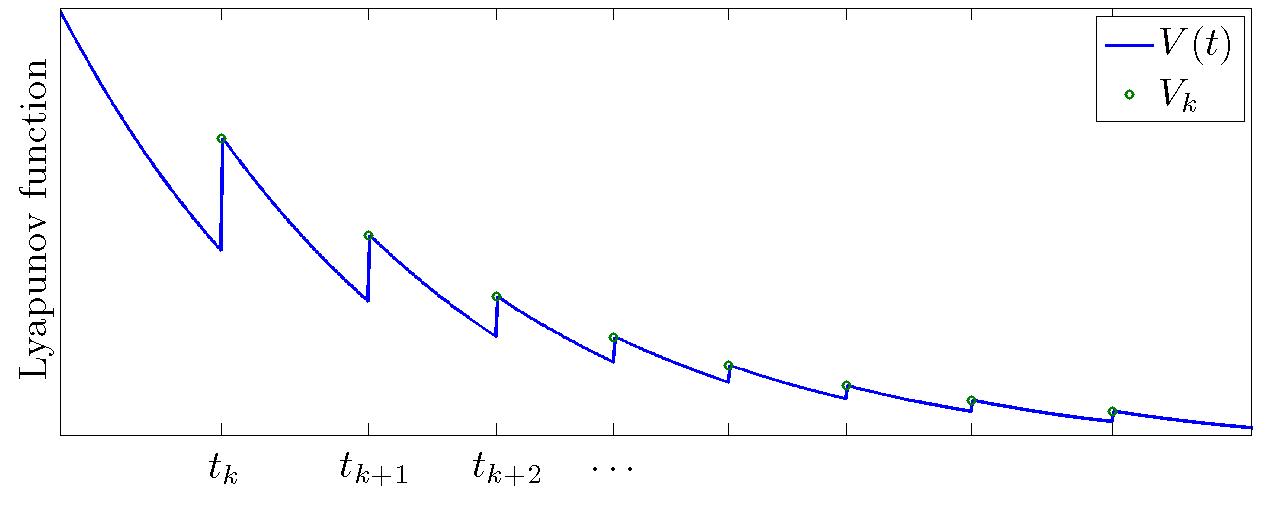}{1}{Time evolution of Lyapunov function $V$.}

\subsection{Boundedness of $L$ and
proof of Theorem \ref{th:main_result}} \label{sec:L_boundedness}
\label{Sec_ProofTh}
Although the construction of the updated law for the high-gain parameter \eqref{eq:high_gain_update_law_start}-\eqref{eq:high_gain_update_law_end} follows the idea developed in \cite{Andrieu2014}, the study of the behavior of the high-gain parameter is more involved. 
Indeed, in the context of observer design of \cite{Andrieu2014}, the nonlinear function $c$ was assumed to be essentially bounded while in the present work, $c$ is depending on $x_1$.
This implies that the interconnection structure between state and high-gain dynamics must be further investigated.

\medskip
\noindent\textbf{Proof of Theorem \ref{th:main_result}:}
Assume $a_1$, $a_2$, $a_3$ and $\alpha^*$ meet the conditions of
Proposition \ref{pr:existenceXk} and Proposition \ref{pr:decreaseVk}.
Consider solutions $(\hat x(\cdot), e(\cdot), L(\cdot),M(\cdot))$
for system \eqref{eq:problem_statement_sys} with the event-triggered output feedback  \eqref{eq:control_law_start}-\eqref{eq:control_law_end} with initial condition $\hat x(0)$ in $\RR^n$, $e(0)$ in $\RR^n$, $L(0))\geq 1$ and $M(0)=1$.
With Proposition \ref{pr:existenceXk} the sequence $(t_k, \xo, \e, \L)_{k\in\NN}$ is well defined.

The existence of a strictly positive dwell time is obtained from the following proposition.
\begin{proposition}\label{Prop_BoundL}
There exists a positive real number $L_{\max{}}$ and class $\mathcal K$ function $\gamma$ and a non decreasing function in both argument $\rho$ such that
\begin{multline}\label{eq:induction_H}
\tilde L_{k+1}
\leq \left(1-\frac{a_1\alpha}{2}\right)\tilde L_k + \gamma(V_k),
\quad \forall k\in \NN,
\end{multline}
where $\gamma(s)=0$ for all $s$ in $[0,1]$ with
$$
\tilde L_k=\max\{L_{k}-L_{\max},0\},
$$
and for all $t$ on the time existence of the solution, we have
\begin{equation}\label{eq_BoundFinalL}
1\leq L(t) \leq \rho(\tilde L_0,V_0)\ .
\end{equation}
\end{proposition}
The proof of this proposition is given in Appendix \ref{sec_ProofProp_BoundL}.

With this proposition in hand, note that it yields for all $k$ in $\NN$,
$\delta_k \ge \frac{\alpha}{\rho(\tilde L_0,V_0)} >0$. 
Consequently, there is a dwell time and the solution are complete (i.e. $\sum_k \delta_k = +\infty$).
Moreover, for all $k$ in $\NN$,
$\frac{L_k}{L^-_{k+1}}\geq \frac{1 }{\rho(\tilde L_0,V_0)}$.
Consequently,
inequality (\ref{eq:LyapDiscrete}) becomes
\[
V_{k+1}
\leq (1-\sigma(\tilde L_0,V_0))
    V_k,
\]
where $\sigma(\tilde L_0,V_0) = \frac{\alpha N(\alpha)}{\rho(L_0,V_0)^{2(n-1+b)}}$ is a decreasing function of both arguments.
This gives $V_k\leq (1-\sigma(\tilde L_0,V_0))^kV_0$, for all $ k$ in $\NN$.
With, \eqref{eq:induction_H},  it yields $\tilde L_k \leq \beta_L(\tilde L_0+V_0,k)$
where
\begin{multline}
\beta_L(s,k)=s\left(1-\frac{a_1\alpha}{2}\right)^k \\+ \sum_{j=1}^k\left(1-\frac{a_1\alpha}{2}\right)^j\gamma\left ((1-\sigma(s,s))^{k-j}s\right).
\end{multline}
The function $\beta_L$ is of class $\mathcal K$ in $s$. Moreover, since $\gamma(s)=0$ for $s\leq 1$, this implies that there exists $k^*(s)$ such that
the mapping $k \mapsto\beta_L(s,k)$ is decreasing for all $k\geq k^*(s)$. Moreover, we have $\lim_{k\rightarrow\infty}\beta_L(s,k) = 0$.
On another hand, since $\delta_k\leq \alpha$, it implies that $k\leq \frac{t}{\alpha}$ for all $t$ in $[t_k, t_{k+1})$. 
\begin{align}
\tilde L(t) &\leq \frac{\tilde L_{k+1}}{1-a_1\alpha}, \notag\\
&\leq \frac{\beta_L(\tilde L_0+V_0,k+1)}{1-a_1\alpha}\label{eq_Lt}.
\end{align}
Finally, 
with \eqref{eq:InterSamp}, it yields
\begin{equation}\label{eq_Vt1}
V(t)\leq k(\alpha)(1-\sigma(\tilde L_0,V_0))^{\frac{t}{\alpha}}V_0\ .
\end{equation}
With the right hand side of \eqref{eq_BoundFinalL} and the definition of the Lyapunov function $V$, we have
\begin{equation}\label{eq_Vt2}
\frac{p_1 + \mu q_1}{2\rho(\tilde L_0,V_0)^{2(n-1+b)}}\left(\norm{ x(t)}^2 + \norm{ \hat x(t)}^2\right)\\
\leq V(t),
\end{equation}
Moreover, we have also:
\begin{equation}\label{eq_Vt3}
V_0 \leq 2 (p_2 + \mu q_2)\left(\norm{ x(0)}^2 + \norm{ \hat x(0)}^2\right).
\end{equation}
From equations \eqref{eq_Lt}, \eqref{eq_Vt1}, \eqref{eq_Vt2}, \eqref{eq_Vt3} and the properties of the function $\beta_L$, it yields readily that there exists a class $\mathcal{KL}$ function $\beta$ such that  inequality \eqref{eq_MainResult} holds.

\section{Illustrative example}
\label{sec_IllustrativeExample}

We apply our approach to the following uncertain third-order system proposed in \cite{Krishnamurthy2004}
\begin{systn} \label{eq:illustrative_example_sys_Krishna}
	\dot{x}_1 &= x_2\\
	\dot{x}_2 &= x_3 \\
	\dot{x}_3 &= \theta x_1^2 x_3 + u,
\end{systn}
where $\theta$ is a constant parameter which only a magnitude bound $\theta_{\max}$ is known. The stabilization of this problem is not trivial even in the case of a continuous-in-time controller. The difficulties come from the nonlinear term $x_1^2 x_3$ that makes $x_3$ dynamics not globally Lipschitz, and from the uncertainty on $\theta$ value, preventing the use of a feedback to cancel the nonlinearity.

However,
system \eqref{eq:illustrative_example_sys_Krishna} belongs to the class of systems \eqref{eq:problem_statement_sys}
and the Assumption \ref{hyp:incremental_bound} is satisfied with $c(x_1) = \theta_{\max} x_1^2$.
Hence,
by Theorem \ref{th:main_result},
an event-triggered output feedback controller \eqref{eq:control_law_start}-\eqref{eq:control_law_end} can be constructed.
Simulation were conducted with gain matrices $K_o$  and $K_c$ and coefficient $\alpha$ selected as
$K_o=\begin{bmatrix}-8 & -12 & -16\end{bmatrix}',\quad K_c= \begin{bmatrix}-15 & -75 &-125 \end{bmatrix}, \quad \alpha=0.1$
to stabilize the linear part of the system \eqref{eq:illustrative_example_sys_Krishna}.\\
Parameters $a_1$, $a_2$ and $a_3$ have then been selected through a trial and error procedure as follows: $$a_1=1, \quad a_2=.5, \quad a_3=.5.$$

Simulation results are given in \FIG{fig:outputFeed_Krishna_example_x_u} and \FIG{fig:outputFeed_Krishna_example_V_delta_L}. 
The evolution of the control and state trajectories are displayed in \FIG{fig:outputFeed_Krishna_example_x_u} for a particular initial condition. 
The corresponding evolution of the Lyapunov function $V$ and the high-gain $L$ are shown in \FIG{fig:outputFeed_Krishna_example_x_u}. We can see how the inter-execution times $\delta_k$ adapts to the nonlinearity. Interestingly, it allows a significant increase of $\delta_k$ when the state is close to the origin: $L(t)$ then goes to $1$ and consequently $\delta_k$ increases toward $\alpha$ value (that was selected as $\alpha=0.1$ in this simulation).

\Figure{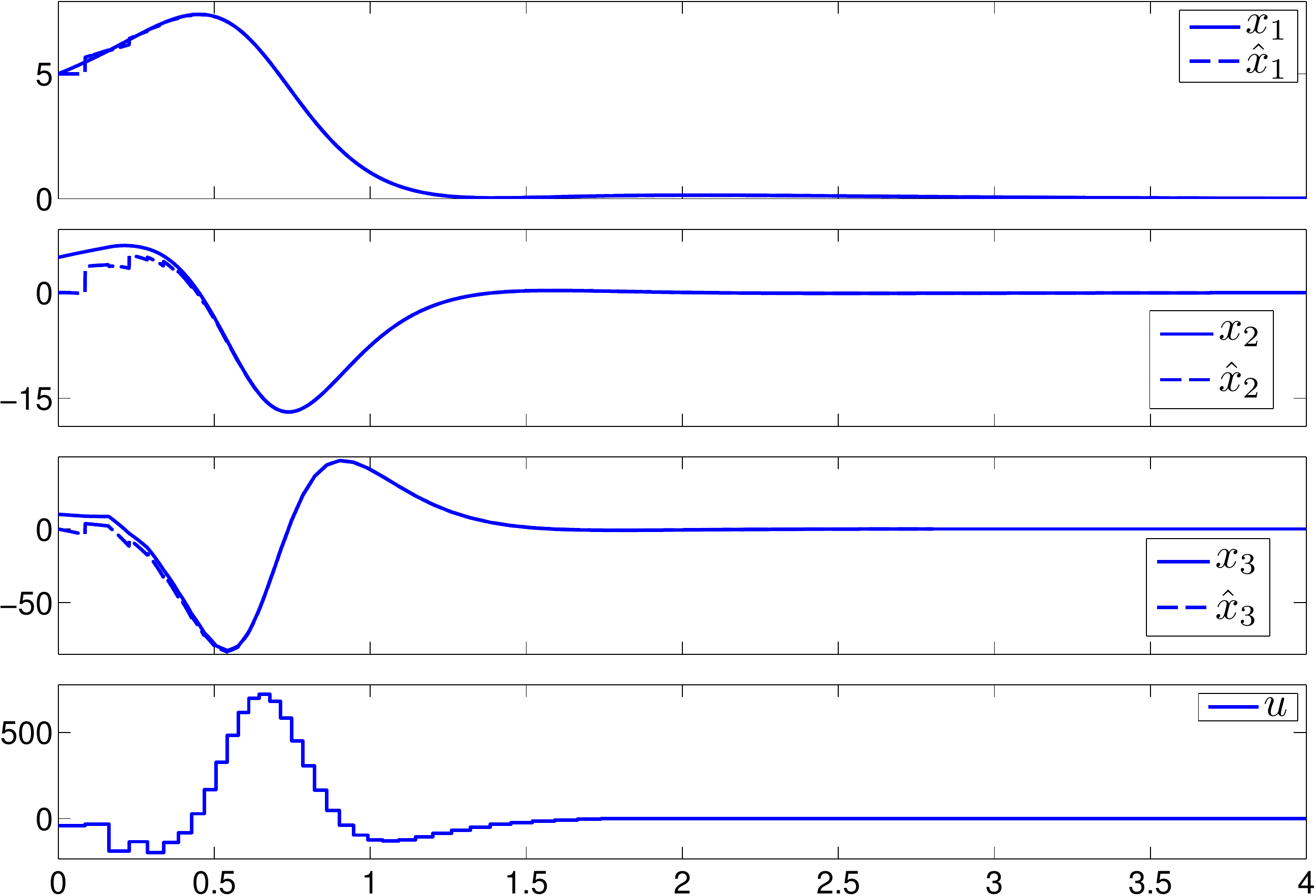}{1}{Control signal and state trajectories of \eqref{eq:illustrative_example_sys_Krishna} with $\left(x_1, x_2, x_3\right) = (5,5,10)$ and $\left(\hat x_1, \hat x_2, \hat x_3\right) = (5,0,0)$ as initial conditions.}
\Figure{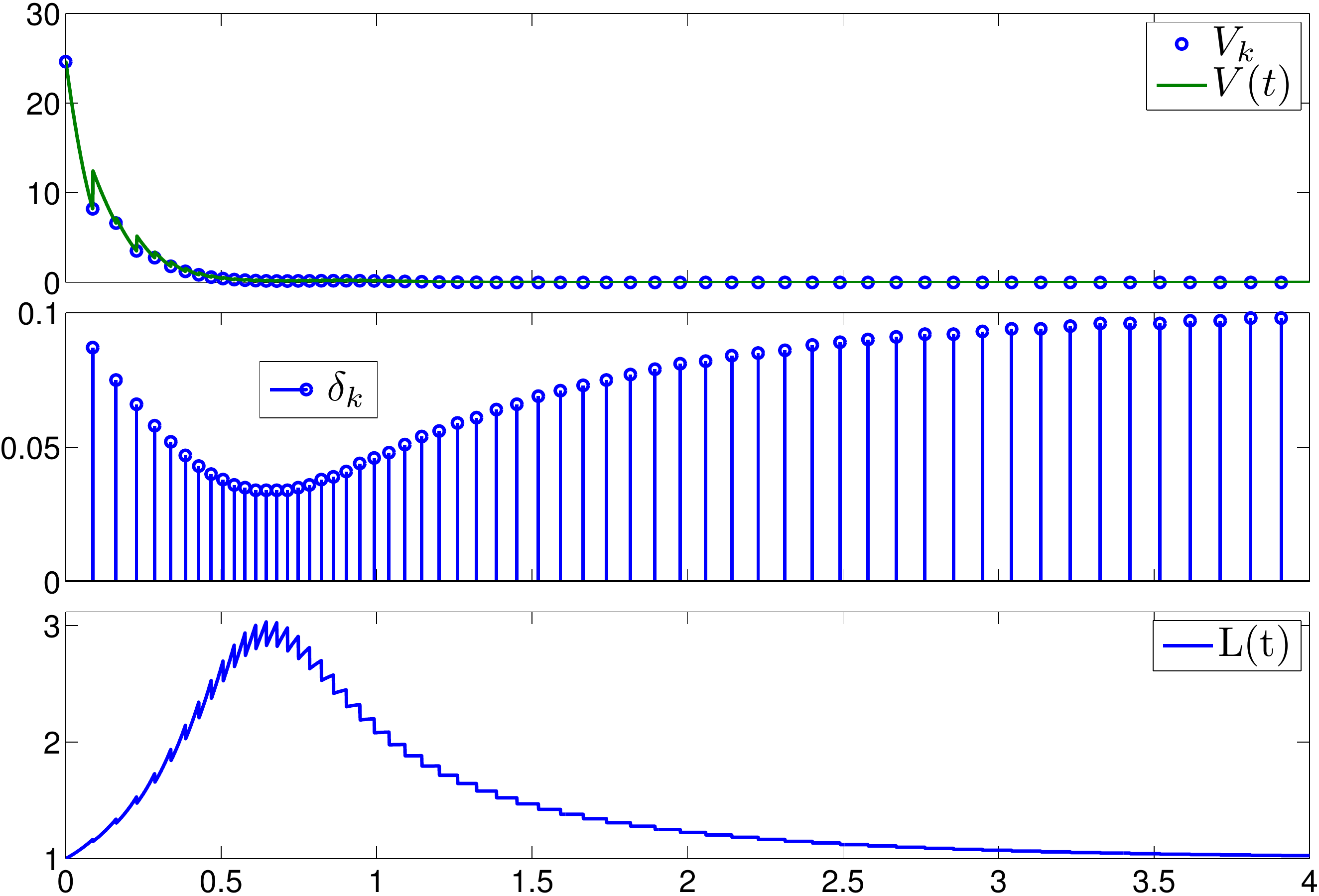}{1}{Simulation results}

\section{Conclusion}
In conclusion, we have presented a new event triggered output feedback for a class of nonlinear systems. 
The triggered mechanism depends on an additional dynamics. This additional dynamics is employed to modify the output feedback following a high-gain paradigm.
The stabilization of the origin of the system is demonstrated and the interest of our approach is illustrated on an example.

\appendix

\section{Proofs of Lemmas}
\subsection{Proof of Lemma \ref{lemm:trigg_linear}} \label{sec:proof_lemma_linear}
The matrix  $(A+BK_c)$ being Hurwitz, let $P$ be a symmetric positive definite matrix such that
\begin{align} \label{eq:lin_K_ineq}
		P(A+BK_c) + (A+BK_c)'P \leq -I,\\
		p_1 I\leq P \leq p_2 I, \notag
\end{align}
with $p_{1},p_{2}$ positive real numbers.
Likewise, let $Q$ be a symmetric positive definite matrix such that
\begin{align}  \label{eq:lin_Ko_ineq}
		Q(A+K_oC) + (A+K_oC)'Q \leq -I,\\
		q_{1} I \leq Q \leq q_{2} I, \notag
\end{align}
with $q_{1},q_{2}$ positive real numbers.

%
%
In order to prove that the origin of the discrete time system \eqref{eq_Discrete} is GAS, we consider the Lyapunov function 
\begin{align}
	V(e,\hat x) = \hat x' P \hat x + \mu e' Q e,
\end{align}
where $\mu$ is a positive real number that will be selected later on.
From \eqref{eq_Discrete}, it comes
\begin{align}
	\eNext'Q  \eNext = e_k' \Fo (\delta)' Q \Fo (\delta) e_k.
\end{align}
Given $v$ in $S^{n-1}=\{v\in\RR^n \mid \|v\|=1 \}$, consider the function
\begin{equation*}
	\nu(\delta,v) = v'\Fo (\delta)'Q\Fo (\delta)v. 
\end{equation*}
We have
\begin{align*}
	\nu(0,v) &= v'Qv, \\
	\dd{\nu}{\delta}(0,v) &= v'[Q(A+K_oC) + (A+K_oC)'Q]v .
\end{align*} 
So, using the inequalities in \eqref{eq:lin_Ko_ineq} , we get
\begin{align}
	\dd{\nu}{\delta}(0,v) &\leq - \dfrac{1}{q_{2}} v'Qv. \label{eq:trigg_linear_dvdtau_ineq}
\end{align} 
Now, we can write
\begin{align*}
	\nu(\delta,v) &= v'Qv + \delta \dd{\nu}{\delta}(0,v) + \rho(\delta,v),
\end{align*}
with $\lim_{\delta \to 0}\frac{\rho(\delta,v)}{\delta} =0$. This equality together with \eqref{eq:trigg_linear_dvdtau_ineq} imply that 
\begin{align*}
	\nu(\delta,v) &\leq (1- \dfrac{\delta}{q_2})v'Qv + \rho(\delta,v).
\end{align*}
The vector $v$ being in a compact set and the function $\rho$ being continuous, there exists $\deltaMaxo$ such that for all $\delta$ in $[0; \deltaMaxo)$ we have $\rho(\delta,v) \leq \frac{\delta}{2 q_2} v'Qv $ for all $v$. This gives
\begin{equation*}
	\nu(\delta,v) \leq \left(1-\frac{\delta}{2 q_2}\right)v'Qv , \qquad \forall~\delta\in[0,\deltaMaxo),
\forall~v\in S^{n-1}. 
\end{equation*}
This property being true for every $v$ in $S^{n-1}$, we have
\begin{align*}
	\Fo (\delta)'Q \Fo(\delta) &\leq \left(1-\frac{\delta }{2  q_2}\right) Q ,
\end{align*}
and there exists $\deltaMaxo$ such that for all $\delta$ in $[0; \deltaMaxo)$ we have 
\begin{align}
\eNext'Q  \eNext 
	\leq \left(1- \dfrac{\delta}{2q_2} \right) e_k' Q  e_k. \label{eq:lin_eNextPo}
\end{align}
Similarly,  we have
\begin{multline}
\xNexto'P \xNexto 
	= \xo' \Fc(\delta)' P \Fc(\delta) \xo + e_k' F_{oc}(\delta)' P F_{oc}(\delta) e_k \notag \\
		 + 2 \xo' \Fc(\delta)' P F_{oc}(\delta) e_k,  \label{eq:lin_xNextoP}
\end{multline}
where $F_{oc}(\delta) = \delta K_o C \exp(A \delta)$.
Notice that $\Fc(0)=I$ and $\dd{\Fc}{\delta}(0)=A+BK_c$. Hence, it implies the existence of a $\deltaMaxc$ such that for all $\delta$ in $[0,\deltaMaxc)$, we have
\begin{align}
\xo' \Fc(\delta)' P \Fc(\delta) \xo 
	\leq & \xo'  P  \xo - \dfrac{\delta}{2p_{2}}\xo' P \xo.
\end{align}
Previous inequality with \eqref{eq:lin_eNextPo} and \eqref{eq:lin_xNextoP} yields
\begin{align*}
V_{k+1} & - V_k  \\
	=&\mu \eNext' Q  \eNext - \mu e_k' Q  e_k +  \xNexto ' P \xNexto  -  \xo ' P \xo \notag \\
	\leq &  - \mu\dfrac{\delta}{2q_{2}}e_k' Q  e_k -  \dfrac{\delta}{2p_{2}}\xo' P \xo +  e_k' F_{oc}(\delta)' P F_{oc}(\delta) e_k \notag\\
	& \pushright{ + 2  \xo' \Fc(\delta)' P F_{oc}(\delta) e_k} \notag \\
	\leq &  - \mu\dfrac{\delta q_{1}}{2q_{2}}\norm{e_k}^2 -  \dfrac{\delta p_1}{2p_{2}}\norm{\xo}^2 +  \norm{F_{oc}(\delta)}^2\norm{P} \norm{e_k}^2 \notag \\
	  & \pushright{ + 2  \norm{\Fc(\delta)}\norm{F_{oc}(\delta)}\norm{P} \norm{\xo} \norm{e_k}.} 
\end{align*}
Using Young's inequality,
the preceding inequality becomes
$$
V_{k+1} - V_k  \leq \left(- \mu\dfrac{\delta q_{1}}{2q_{2}} +  N(\delta) \right) \norm{e_k}^2 -  \dfrac{\delta p_1}{4p_{2}} \norm{\xo}^2
$$
where
$$N(\delta) = \norm{F_{ex}(\delta)}^2\norm{P}	+  \norm{F_x(\delta)}^2\norm{F_{ex}(\delta)}^2\norm{P}^2 \dfrac{4p_{2}}{\delta p_1}.$$
Then, choosing $\mu$ as
\begin{align*}
	\mu \geq \dfrac{2q_{2}N(\delta)}{\delta q_{1}},
\end{align*}
ensures the decrease of $V$ for all $\delta$ in $[0,\deltaMax)$, with $\deltaMax = \max\{\deltaMaxc, \deltaMaxo\}$.

\subsection{Proof of Lemma \ref{lemm:state_Feed_To2}}
\label{sec:Proof_Lemm_To2}
The proof of Lemma \ref{lemm:state_Feed_To2} uses  \cite[Lemma 6, p112]{Andrieu2014}.
\begin{lemma}[ \cite{Andrieu2014}]\label{lemm:state_feed_PsiPPsi}
The matrix $Q$ and $P$ satisfy the following property for all $a_1$ and $\alpha$ such that $a_1\alpha<1$
$$
\Psi Q \Psi\le \psi_0(\alpha)Q\psi_0(\alpha)\ , \Psi P \Psi\le \psi_0(\alpha)P\psi_0(\alpha),
$$
where $\Psi=\Snext \inv{\SnextM}$ and
\[
\psi_0(\alpha)=\diag\left(\frac{1}{(1-a_1\alpha)^b}, \dots, \frac{1}{(1-a_1\alpha)^{n+b-1}}\right).
\]
\end{lemma}

\vspace{0.5em}
To prove Lemma \ref{lemm:state_Feed_To2}, we first analyse the term $R_{o}$.
Following what has been done in \eqref{eq_SF}, it yields
\begin{equation}
\norm{\SnextM f(x(t_k+s))}^2\leq n^2 c(t_k+s)^2 \norm{\SnextM x(t_k+s)}^2 . \label{eq:temp}
\end{equation} 
From the previous inequality, we get
\begin{multline}
\norm{R_o } \leq  \norm{I+ \alpha K_o C}  \exp(\norm{A} \alpha) \\
 \times \int_{0}^{\delta_k}n c(t_k+s) \norm{\SnextM x(t_k+s)} ds. \label{eq:outputFeed_normRo}
\end{multline}
On the other hand, we have for all $s$ in $[0,\delta_k)$
\begin{align*}
	\SnextM \dot{x} (t_k+s)
		=& \SnextM \big(Ax(t_k+s) + B K_c (L_k)^{n+1}\mathcal{L}_k  \xo \\
		&	\pushright{ + f(x(t_k+s)) \big)} \\
		=& \LnextM A \SnextM x(t_k+s)  \\
		&	+ \LnextM  B K_c \Omega \SnextM \xo + \SnextM f(x(t_k+s)).
\end{align*}
where
\begin{align*}
\Omega
		=& (\LnextM)^{-n-1}(\L)^{n+1}  \Lmat \inv{\mathcal L_{k+1}^-} \\
		=& \diag\left\{\left(\frac{L_k}{L_{k+1}^-}\right)^n, \left(\frac{L_k}{L_{k+1}^-}\right)^{n-1},\dots,\frac{L_k}{L_{k+1}^-}\right\}
\end{align*}
Note that since $L_{k+1}^-\geq L_k$, it yields $\norm{\Omega}\leq 1$. Hence,denoting by $w(s)$ the expression $\SnextM x(t_k+s)$, this gives
\begin{align*}
\dfrac{d}{ds} \norm{w(s)} 
	=& \dfrac{\langle \dot{w}(s),w(s)\rangle}{\norm{w(s)}} \\
	\leq& \left(\LnextM \norm{A}+n c(t_k+s)\right) \norm{ w(s)}\\
		&\pushright{	 +\LnextM \norm{B K_c}\norm{\SnextM \xo}}
\end{align*}
Hence, by integrating preceding inequality, it yields
\begin{multline*}
	\norm{w(s)} 
 \leq \int_{0}^{s} (\LnextM \|A\| + nc(t_k+r)) \|w(r)\| dr \\
	 + \LnextM \norm{B K_c} \norm{\SnextM \xo}  s + \norm{w(0)} .
\end{multline*}
Since $(\LnextM \|A\| + nc(t_k+s))$ is a continuous non-negative function and $\LnextM \| B K\|\| \SnextM x_k \| s + \norm{w(0)}$ is non-decreasing, applying a variant of the Gronwall-Bellman inequality \cite[Theorem 1.3.1]{ames1997inequalities}, it comes
\begin{multline}\label{eq:outputFeed_Gronwall_fin}
	\norm{\SnextM x(t_k+s)}
	 \leq (\LnextM \norm{B K_c} \norm{\SnextM \xo}  s + \norm{\SnextM x_k}) \\
		 \times	 \exp\left(\LnextM \|A\| s\right)\exp\left(\int_{0}^{s}( nc(t_k+r)dr\right),	
\end{multline}
Consequently, according to \eqref{eq:outputFeed_normRo}, we get
\begin{align*}
	 \norm{\Ro }& \\	
	\leq& \norm{I+ \alpha K_o C} \exp(2\norm{A}\alpha) \notag \\
	& \times \int_{0}^{\delta_k}n c(t_k+s) (\LnextM\norm{ B K_c} \norm{\SnextM \xo}s + \norm{\SnextM x_k}) \\
		&\pushright{	\times \exp\left(\int_{0}^{s}nc(t_k+r)dr\right)ds } \\
	\leq& \norm{I+ \alpha K_o C} \exp(2\norm{A}\alpha) \left(\alpha \norm{ B K_c} \norm{\SnextM \xo} +\norm{\SnextM x_k} \right) \\
	&\pushright{	 \times \int_{0}^{\delta_k}n c(t_k+s) \exp\left(\int_{0}^{s}nc(t_k+r)dr\right) ds }\\
	=& \norm{I+ \alpha K_o C} \exp(2\norm{A}\alpha) \left(\alpha \norm{ B K_c} \norm{\SnextM \xo} +\norm{\SnextM x_k} \right) \\
		&\pushright{	\times  \left[ \exp\left(\int_{0}^{s}nc(t_k+r)dr\right) \right]_{s=0}^{s=\delta_k}}\\
	=& \norm{I+ \alpha K_o C} \exp(2\norm{A}\alpha) \left(\alpha \norm{ B K_c} \norm{\SnextM \xo} +\norm{\SnextM x_k} \right) \\
		&\pushright{	\times  \left[\exp\left(\int_{0}^{\delta_k}nc(t_k+r)dr\right)-1\right] .}
\end{align*}
Hence, employing $e_k = \hat x_k -x_k$  it yields,
\begin{multline*}
\norm{\Ro }\leq \left[M_{2}(\alpha)\norm{\SnextM \xo} + M_{1}(\alpha)\norm{\SnextM e_k}\right]  \\
\times\left[e^\beta
-1\right] .
\end{multline*}
where
\begin{align*}
M_{1}(\alpha) &= \norm{I+ \alpha K_o C} \exp(2\norm{A}\alpha),\\
M_{2}(\alpha) &= M_{1}(\alpha)\left(\alpha \norm{ B K_c}  +1 \right).
\end{align*}
Hence, employing Lemma \ref{lemm:state_feed_PsiPPsi} this gives
\begin{multline*}
\norm{\Ro'\Psi P \Psi\Ro}\leq 
\left[M_{4}(\alpha)\norm{\SnextM \xo}^2 + M_{3}(\alpha)\norm{\SnextM e_k}^2\right] \\
\times [e^{\beta}-1]^2
\end{multline*}
where
\begin{align*}
M_{3}(\alpha) &= \dfrac{2\norm{Q}}{(1-a_1\alpha)^{2(n-b+1) }}M_{e,1}(\alpha)^2,\\
M_{4}(\alpha) &= \dfrac{2\norm{Q}}{(1-a_1\alpha)^{2(n-b+1) }}M_{\hat x,1}(\alpha)^2.
\end{align*}
Moreover,
\begin{multline*}
\norm{2 e_k' \SnextM \Fo (\alpha)' \Psi Q \Psi \Ro}\leq\\
\left[M_{6}(\alpha)\norm{\SnextM \xo}^2 + M_{5}(\alpha)\norm{\SnextM e_k}^2\right] \\
\times [e^{\beta}-1],
\end{multline*}
where
\begin{align*}
M_{5}(\alpha) &= \dfrac{2\norm{Q}\norm{\Fo}}{(1-a_1\alpha)^{2(n-b+1) }}M_{e,1}(\alpha),\\
M_{6}(\alpha) &= \dfrac{\norm{Q}\norm{\Fo}}{(1-a_1\alpha)^{2(n-b+1) }}M_{\hat x,1}(\alpha).
\end{align*}
Noticing that 
\begin{equation}\label{eq_ebeta}
0\leq (e^\beta-1)^2\leq e^{2\beta} - 1\ ,\ 0\leq (e^\beta-1)\leq e^{2\beta} - 1,
\end{equation}
the result follows with
$$
N_{o,e}(\alpha) = M_{3}(\alpha) + M_{5}(\alpha)\ ,\
N_{o,\hat x}(\alpha) = M_{4}(\alpha) + M_{6}(\alpha).
$$

\subsection{Proof of Lemma \ref{lemm:state_Feed_Tc2}}
\label{sec:Proof_Lemm_T2}
The first part of the proof is devoted to upper-bound the term $\norm{\Rc} = \alpha \norm{  K_o C \SnextM e_{k+1}^-}$.
From the algebraic equality given in \eqref{eq:linear_algebraic_prop} and the expression of $e_{k+1}^-$ given in \eqref{eq:outputFeed_eNextM}, it yields
\begin{multline*} 
\norm{\Rc} \leq M_7(\alpha) \Big[ \norm{\SnextM e_k}\\
+ \int_{0}^{\delta_k}n c(t_k+s) \norm{\SnextM x(t_k+s)} ds\Big],
\end{multline*}
where $M_7(\alpha) = \alpha \norm{ K_o C}\exp(|A|\alpha)$.
Consequently, according to \eqref{eq:temp} and \eqref{eq:outputFeed_Gronwall_fin}, we get
\begin{align*}
	\norm{\Rc} 
	\leq& M_7(\alpha)  \Big[ \norm{\SnextM e_k}+	\exp(\norm{A}\alpha) \int_{0}^{\delta_k}n c(t_k+s)\\
	& \times  (\LnextM\norm{ B K_c} \norm{\SnextM \xo}s + \norm{\SnextM x_k})\Big] \\
	& \times \exp\left(\int_{0}^{s}(nc(t_k+r)dr\right)ds, \\
	\leq& M_7(\alpha) \Big[ \norm{\SnextM e_k}+	\exp(\norm{A}\alpha) \\
	& \times	 \int_{0}^{\delta_k}n c(t_k+s)\exp\left(\int_{0}^{s}(nc(t_k+r)dr\right)ds   \\
	& \times  (\alpha\norm{ B K_c} \norm{\SnextM \xo} + \norm{\SnextM x_k})\Big] ,\\
	=&  M_7(\alpha)\Big[ \norm{\SnextM e_k}  +	(e^\beta-1) \exp(\norm{A}\alpha)  \\
	& \times (\alpha\norm{ B K_c} \norm{\SnextM \xo} + \norm{\SnextM x_k})\Big] .
\end{align*}
Hence, employing $e_k = \hat x_k -x_k$  it yields,
\begin{multline*}
\norm{\Rc }\leq \left[M_8(\alpha)\norm{\SnextM \xo} + M_9(\alpha)\norm{\SnextM e_k}\right]\left[e^\beta
-1\right]  \\
+ M_7(\alpha)\norm{\SnextM e_k}  .
\end{multline*}
where
\begin{align*}
M_8(\alpha) &= M_7(\alpha) (\alpha\norm{ B K_c}+1) \exp(\norm{A}\alpha) ,\\
M_9(\alpha) &= M_7(\alpha)  \exp(\norm{A}\alpha).
\end{align*}
Hence, employing Lemma \ref{lemm:state_feed_PsiPPsi} and \eqref{eq_ebeta} this gives
\begin{multline}\label{Ineq_1}
\norm{\Rc'\Psi Q \Psi\Rc}\leq 
\left[M_{10}(\alpha)\norm{\SnextM \xo}^2 + M_{11}(\alpha)\norm{\SnextM e_k}^2\right] \\
\times [e^{2\beta}-1] + M_{12}(\alpha)\norm{\SnextM e_k}^2
\end{multline}
where
\begin{align*}
M_{10}(\alpha) &= \dfrac{3\norm{P}M_{8}(\alpha)^2}{(1-a_1\alpha)^{2(n-b+1) }},\\
M_{11}(\alpha) &= \dfrac{\norm{P}\left[2M_{9}(\alpha)^2+M_7(\alpha)^2 + 2 M_9(\alpha)M_7(\alpha)\right]}{(1-a_1\alpha)^{2(n-b+1) }},\\
M_{12}(\alpha) &= \dfrac{\norm{P}M_{7}(\alpha)^2}{(1-a_1\alpha)^{2(n-b+1) }}.
\end{align*}

On another hand, with the algebraic equality given in \eqref{eq:linear_algebraic_prop},  we have\\[0.5em]
$\displaystyle
\SnextM  \inv{\S}\Fc(\alpha_k)\hat X_k =  $\hfill \refstepcounter{equation} ($\theequation$)\label{eq_NewProof}\\[0.5em]
\null\hfill$\displaystyle
\left[\exp(A\alpha)  + \int_0^\alpha\exp(A(\alpha-s))dsBK_c \Lambda \right] \SnextM \hat x_k,$\\[0.5em]
where
$ \Lambda =
\left(\frac{L_k}{L_{k+1}^-}\right)^{n+1} \S  \inv{\SnextM}$.
Note that $L_{k+1}^-\geq L_k$. Hence,  $\norm{\Lambda} \leq 1$ and  we have
\begin{multline*}
\norm{\SnextM  \inv{\S}\Fc(\alpha_k)\hat X_k} \leq  \\\left[\exp(\norm{A}\alpha)(1  + \norm{BK_c})\right] \norm{\SnextM \hat x_k}.
\end{multline*}
Hence, employing Lemma \ref{lemm:state_feed_PsiPPsi}, this gives
\begin{multline}\label{Ineq_2}
 2 \Rc' \Psi P\mathcal S_{k+1} \inv{\mathcal S_{k} }\Fc (\alpha_k) \hat X_k \leq \\ \left[M_{13}(\alpha)\norm{\SnextM \xo}^2 + M_{14}(\alpha)\norm{\SnextM e_k}^2\right] [e^{2\beta}-1] \\
+ M_{15}(\alpha)\norm{\SnextM e_k}\sqrt{U_k}
\end{multline}
where
\begin{align*}
M_{13}(\alpha) &= \dfrac{\norm{P}(M_{8}(\alpha)+\frac{1}{2})\left[\exp(\norm{A}\alpha)(1  + \norm{BK_c})\right] }{(1-a_1\alpha)^{2(n-b+1) }},\\
M_{14}(\alpha) &= \dfrac{\norm{P} M_{9}(\alpha)}{2(1-a_1\alpha)^{2(n-b+1) }},\\
M_{15}(\alpha) &= \dfrac{\norm{P}M_{7}(\alpha)\left[\exp(\norm{A}\alpha)(1  + \norm{BK_c})\right]}{\sqrt{\norm{P}}(1-a_1\alpha)^{2(n-b+1) }}.
\end{align*}
and where we have used $\norm{\SnextM \hat x_k} = \norm{\SnextM \inv{\S}\hat X_k}\leq \sqrt{\frac{U_k}{\norm{P}}}$
Finally note that 
\begin{multline*}
M_{15}(\alpha)\norm{\SnextM e_k}\sqrt{U_k}\leq \frac{1}{2}\left(\frac{\alpha}{p_2}\right)^2U_k \\+ \frac{1}{2}\frac{ M_{15}(\alpha) p_2^2}{\alpha^2}\norm{\SnextM e_k}.
\end{multline*}
Hence the result follows from the former inequality in combination with inequalities  \eqref{Ineq_1} and  \eqref{Ineq_2}.

\subsection{Proof of Proposition \ref{Prop_BoundL}}
\label{sec_ProofProp_BoundL}
\begin{proof}
Inequality (\ref{eq:LyapDiscrete}) of Proposition \ref{pr:decreaseVk} implies that
$(V_k)_{k\in\NN}$ is a nonincreasing sequence.
Consequently, being nonnegative, $(V_k)_{k\in\NN}$ is bounded.
One infers, using inequality (\ref{eq:InterSamp}), that $V(t)$ is bounded.
Hence, by the left parts in inequalities \eqref{eq:bound_P}-\eqref{eq:bound_Po},
we get that, on the time $T_x$ ($=\sum\delta_k$) of existence of the solution,
$\hat X(t)$ and $E(t)$ (and consequently so are $\frac{\hat x_1(t)}{L(t)^b}=\hat X_1(t)$ and $\frac{e_1(t)}{L(t)^b} = E_1(t)$)
are bounded. Then we get that $\frac{ x_1(t)}{L(t)^b}$ is bounded since we have $| x_1(t)| \leq |\hat x_1(t)| + |e_1(t)|$.

Summing up, there exists a class $\mathcal K$ function $\dr_1$  such that
\begin{equation}\label{eq:X1_bounded}
\frac{|x_1(t)|}{L(t)^b}
\le \dr_1(V_k)\le \dr_1(V_0),\quad \forall (t,k)\in [t_k,T^*).
\end{equation}
With this result in hand, let us analyze the high-gain dynamics.
According to equations (\ref{eq:high_gain_update_law_start}) and (\ref{eq:high_gain_update_law_M}),
we have,
for all $k$ and all $t$ in $[t_k, t_{k+1})$,
$\dot L(t) = \frac{a_2}{a_3} L(t) \dot M(t)$,
which implies that
for all $t$ in $[t_k, t_{k+1})$
\begin{align}
L(t) 
	=& \exp\left (\frac{a_2}{a_3} \int_{t_k}^t \dot M(s)ds\right )L_k, \notag \\
	=& \exp\left(\frac{a_2}{a_3} M(t)-\frac{a_2}{a_3}\right)L_k . 
	\label{eq:boundedness_Lcontinu}
\end{align}
Consequently, from \eqref{eq:high_gain_update_law_Lnext} and \eqref{eq:control_law_end}
\begin{align}
\Lnext = \exp\left(\frac{a_2}{a_3}(M\NextM-1)\right)\L (1-a_1\alpha)+a_1\alpha, \label{eq:boundedness_Lnext}
\end{align}
and $\delta_k$ satisfies $$\exp\left(\frac{a_2}{a_3}(M\NextM-1)\right)\delta_k L_k=\alpha.$$
Since $M\NextM\geq 1$, $a_2\geq 0$ and $a_3\geq 0$ the previous equality implies
\begin{equation}\label{2}
\delta_k L_k\leq \alpha .
\end{equation}
Moreover, we have
\begin{align}
\dot M(t) &= a_3 M(t)c(x_1(t)) \notag\\
&= a_3 M(t)(c_0 + c_1 |x_1|^q) \notag\\
&\leq a_3 M(t) ( c_0 + c_1 \dr_1(V_k)^q L(t)^{bq}) \tag{by (\ref{eq:X1_bounded})} \\
&\leq a_3 (c_0 + c_1 \dr_1(V_k)^q)  M(t) L(t)^{bq} \tag{since $L(t)\geq 1$} \\
&\leq \dr_2(V_k)  M(t) \exp\left(\frac{a_2}{a_3}bq(M(t)-1)\right)L_k^{{bq}}, \tag{by (\ref{eq:boundedness_Lcontinu})}
\end{align}
where $\dr_2(V_k)=a_3(c_0 + c_1 \dr_1(V_k)^q)$.
Let $\psi(t)$ be the solution to the scalar dynamical system
\[
\dot \psi(t) =  \psi(t) \exp\left({\frac{a_2}{a_3}bq}(\psi(t) -1)\right),\quad
\psi(0)=1.
\]
$\psi(\cdot)$
is defined on $[0,T_\psi)$ where $T_\psi$ is a positive real number possibly equal to $+\infty$.
Note that we have
(see e.g. \cite[Theorem 1.10.1]{Lakshmikantham1969})
that for all $t$ such that $0\leq \dr_2(V_k)(t-t_k)L_k^{{bq}}<T_\psi$
$$
M(t) \leq \psi\left(\dr_2(V_k)(t-t_k)L_k^{{bq}}\right).
$$
Consequently,
for all $k$ such that $\dr_2(V_k)\delta_k L_k^{{bq}}<T_\psi$
\[
M\NextM = M(t_k+\delta_k^-) \leq \psi\left (\dr_2(V_k)\delta_k L_k^{bq}\right ).
\]
From this,
we get employing (\ref{2}) that,
for all $k$ such that $\dr_2(V_k)\alpha L_k^{{{bq}}-1}<T_\psi$
\begin{equation}
1\leq M\NextM\leq \psi\left ( \dr_2(V_k)\alpha L_k^{bq-1}\right ), \label{eq:boundedness_Mnext_bounds}
\end{equation}
and employing \eqref{eq:boundedness_Lnext} that,
for all $k$ such that $\dr_2(V_k)\alpha L_k^{{{bq}}-1}<T_\psi$
\begin{equation}\label{eq_usebq}
L_{k+1} \leq F(L_k) ,
\end{equation}
where
$$
F(L_k)= \exp\left (\psi\left ( \dr_2(V_k) \alpha L_k^{bq-1}\right )-1\right )L_k (1-a_1\alpha)+a_1\alpha.
$$
Note that, since $bq<1$,
\[
\lim_{L\to+\infty}L^{bq-1} = 0
\]
and since moreover, $\psi(0)=1$, we also get
\[
\lim _{L\rightarrow +\infty} \frac{F(L)}{L} = 1-a_1\alpha<1.
\]
Consequently,
there exists an increasing function $\bar L_1$ such that for all $L> \bar L_1(V_k)$
\begin{equation}\label{eq_barL}
\dr_2(V_k)\alpha{L}^{bq-1}<T_\psi,
\quad
F(L) < \left(1-\frac{a_1\alpha}{2}\right) L.
\end{equation}
On the other hand, consider the following nonlinear system with input $\chi$
\begin{equation}\label{eq_LM}
\left\{
\begin{aligned}
\dot L(t)
&= a_2L(t)M(t) \left(c_0 + c_1\chi(t){^q} L(t)^{bq}\right)\\
\dot M(t)
&= a_3 M(t)\left(c_0 + c_1\chi(t){^q} L(t)^{bq}\right),
\end{aligned}\right.
\end{equation}
We assume that the norm of the input signal satisfies the bound 
\begin{equation}\label{eq_BoundInput}
|\chi(\cdot)|\leq \dr_1(v)\ ,
\end{equation}
where $v$ is a given positive real number.
Notice that the couple $(L,M)$ which satisfies equations (\ref{eq:high_gain_update_law_start}) and (\ref{eq:high_gain_update_law_M}) 
between $[t_k,t_{k+1})$
is also a solution of the previous nonlinear system with input $\chi(t) = \frac{x_1(t)}{L(t)^b}$ which satisfies \eqref{eq_BoundInput} with $v = V_k$.
Let $\phi_{s,t}$ denotes the flow of (\ref{eq_LM}) issued from $s$,
i.e., $\phi_{s,t}(a,b)$ is the solution of (\ref{eq_LM}) that takes value $(a,b)$ at $t=s$.
Let $C_1$, $C_2$, be the two compact subsets of $\RR^2$ defined by:
\begin{align*}
C_1 &=
 \{1\le L\leq \bar L_1(v), M=1\},\\
C_2 &= \{1\leq L \leq 2\bar L_1(v), 0\leq M\leq 2\}.
\end{align*}
The set $C_1$ is included in the interior of $C_2$, and we have the following Lemma.
\begin{lemma}\label{lemm_C1C2}
There exists a non increasing function $t_1$ such that for all input function $\chi$ which satisfies the bound \eqref{eq_BoundInput} the following holds.
\begin{equation}\label{eq_t1}
\forall k\in\NN,\quad \forall t\le t_1(v),\quad \phi_{t_k,t_k+t}(C_1)\subset C_2.
\end{equation}
\end{lemma}
The proof of Lemma \ref{lemm_C1C2} is given in Appendix \ref{sec_Prooflemm_C1C2}.
Let
$$
\bar L_2(v):=\max\Big\{2\bar L_1(v), \frac{\alpha}{t_1(v)}\Big\}\ .
$$
Note that $L_k$ satisfies the following property:
\begin{enumerate}
\item If $L_{k} > \bar L_1(V_k)$ then $L_{k+1}\leq \left(1-\frac{a_1\alpha}{2}\right) L_k$;
\item If $L_{k} \leq \bar L_1(V_k)$ then $L_{k+1}\leq \bar L_{2}(V_k)$.
\end{enumerate}
Indeed, we have
\begin{enumerate}
\item  If $L_{k} > \bar L_1(V_k)$.
With (\ref{eq_usebq}) and (\ref{eq_barL}), we get
\begin{align*}
L_{k+1} 
&\leq \left(1-\frac{a_1\alpha}{2}\right)L_k\ .
\end{align*}
\item If $ L_{k} \le \bar L_1(V_k)$
\begin{enumerate}
\item     If $ \delta_{k}\leq t_1(V_k)$.
Because $L_{k+1}^{-}\ge 1$ and $a_1\alpha<1$, (\ref{eq:high_gain_update_law_Lnext}) implies that
\(
L_{k+1}
\le   L_{k+1}^{-}.
\)
It follows, using (\ref{eq_t1}) with $v=V_k$ (note that $(L_{k},M_{k})\in C_1$), that
\begin{multline*}
L_{k+1}
\le   L_{k+1}^{-}
=   L\left( (t_{k}+\delta_{k})^{-} \right)\\
\le 2 \bar L_1(V_k)
\le \bar L_{2}(V_k).
\end{multline*}
\item If $\delta_{k}> t_1(V_k)$.
\(
L_{k+1}
\le   L_{k+1}^{-},
\)
and since, by (\ref{eq:control_law_end}), ${\delta_k}L_{k+1}^{-}={\alpha}$,
it follows that
\[
L_{k+1}
\le   \frac{\alpha}{\delta_k}
\le \frac{\alpha}{t_1(V_k)}
\le \bar L_{2}(V_k).
\]
\end{enumerate}
\end{enumerate}
Note that the previous properties, implies that for all $k$
$$
L_{k+1} \leq  \left(1-\frac{a_1\alpha}{2}\right) L_k+ \bar  L_{2}(V_k)
$$
and the first part of the result (i.e. inequality \eqref{eq:induction_H}) holds with $L_{\max} = \bar L_{2}(1)$ and  $\gamma(V_k)=\max\{\bar L_{2}(V_k)-\bar L_{2}(1),0\}$.

Note that the previous properties 1) and 2) in combination with the fact that the sequence $(V_k)$ is decreasing imply also for all $k$
$$
L_{k} \leq \max\{ L_0, \bar L_{2}(V_0)\}
$$
Moreover, since for all $k$ in $\NN$ and all $t$ in $[t_k, t_{k+1})$
\begin{align}
L(t)
&\leq L_{k+1}^- \tag{since $\dot L(t)\geq 0$}\\
&= \frac{L_{k+1}-a_1\alpha}{1-a_1 \alpha} \tag{by (\ref{eq:high_gain_update_law_Lnext})}\\
&\leq \frac{L_{k+1}}{1-a_1 \alpha} \notag \\
&\leq \frac{\left(\frac{a_1\alpha}{2}\right)^{k+1}\max\{ L_0, \bar L_{2}(V_0)\}}{1-a_1 \alpha},  \label{eq:bound_Lk}\\
&\leq \frac{\max\{ L_0, \bar L_{2}(V_0)\}}{1-a_1 \alpha},  \label{eq:bound_L}
\end{align}
and the result holds with $\rho(L_0,V_0) = \frac{\max\{ \tilde L_0 + \bar L_{2}(1), \bar L_{2}(V_0)\}}{1-a_1 \alpha}$.
\end{proof}

\subsection{Proof of Lemma \ref{lemm_C1C2}}
\label{sec_Prooflemm_C1C2}
Let $dL_{\max}$ and $dM_{\max}$ be the increasing functions 
\begin{align*}
dL_{\max}(v) &= 4a_2\bar L_1(v) (c_0 + c_1\dr_1(v)^q(2\bar L_1(v))^{bq} ),\\
dM_{\max}(v) &= 2a_3 (c_0 + c_1\dr_1(v)^q(2\bar L_1(v))^{bq} ).
\end{align*}
Note that if $(L(t),M(t))$ is in $C_2$ and $\chi(t)$  satisfies the bound \eqref{eq_BoundInput}, we have
\begin{equation}\label{eq_dLmax}
\dot L(t) \leq dL_{\max}(v) \ ,\ \dot M(t) \leq dM_{\max}(v) .
\end{equation}
Let $t_1$ be the function defined by
$$
t_1(v) = \min\left\{
\frac{1}{dL_{\max}(v)}, \frac{1}{dM_{\max}(v)}\right\}.
$$
We show that this function satisfies the properties of Lemma \ref{lemm_C1C2}.
Assume this is not the case. 
Hence, there exists $M(t_k), L(t_k)$ in $C_1$, $\chi$ which  satisfies the bound \eqref{eq_BoundInput} and $t^*\leq t_1(v)$ such that
$(L(t_k+t^*),M(t_k+t^*))\notin C_2$.
Let $s^*$ be the time instant at which the solution leaves $C_2$. More precisely, let
$s^* = \inf\{s, t_k\leq s\leq t_k+t^*, (L(s),M(s))\notin C_2\}$.
Note that $(L(s^*),M(s^*))$ is at the border of $C_2$ and $t_k<s^*<t_k+t_1(v)$.
Moreover, with (\ref{eq_dLmax}), it yields:
$$
M(s^*) \leq 1 + (s^*-t_k) dM_{\max}(v) <  1+ t_1(v) dM_{\max}(v)  \leq  2.
$$
Similarly, we have
$$
L(s^*)  <  L(t_k) + t_1(v) dL_{\max}  \leq  L(t_k)+1 \leq 2\bar L_1(v),
$$
where the last inequality is obtained since $\bar L_1(v)\geq 1$.
This implies that $(L(s^*),M(s^*))$ is not at the border of $C_2$ which contradicts the existence of $t^*$.

\bibliographystyle{plain}
\bibliography{livres,event_trigg}

\if \hal1
\else
\begin{IEEEbiography}[{\includegraphics[width=1in,height=1.25in,clip,keepaspectratio]{BioPhoto/Photo_Peralez.png}}]%
{Johan Peralez}
received the MSc degree and the PhD degree in Systems and Controls in 2011 and 2015, respectively, from Lyon University of Science, France. 
His thesis work, conducted in part within IFP Energies Nouvelles, focused on Rankine cycle aboard vehicles. 
Between march 2015 and february 2016, he was a postdoctoral Researcher in LAGEP laboratory, France, where he was interested in triggered control theory.
\end{IEEEbiography}
\begin{IEEEbiography}[{\includegraphics[width=1in,height=1.25in,clip,keepaspectratio]{BioPhoto/Photo_Andrieu.jpeg}}]%
{Vincent Andrieu}
 graduated in applied mathematics from “INSA de Rouen”, France, in 2001. After working in ONERA (French aerospace research company), he obtained a PhD degree from “Ecole des Mines de Paris” in 2005. In 2006, he had a research appointment at the Control and Power Group, Dept. EEE, Imperial College London. In 2008, he joined the CNRS-LAAS lab in Toulouse, France, as a “CNRS-charg\'e de recherche”. Since 2010, he has been working in LAGEP-CNRS, Universit\'e de Lyon, France. In 2014, he joined the functional analysis group from Bergische Universit\"at Wuppertal in Germany, for two years. His main research interests are in the feedback stabilization of controlled dynamical nonlinear systems and state estimation problems. He is also interested in practical application of these theoretical problems, and especially in the field of aeronautics and chemical engineering.
\end{IEEEbiography}
\begin{IEEEbiography}[{\includegraphics[width=1in,height=1.25in,clip,keepaspectratio]{BioPhoto/Photo_Nadri.png}}]%
{Madiha Nadri}
 received the MSc degree in Automatic Control from universit\'e Claude Bernard of Lyon in 1997. In 2001, she received the PhD degree from the same university. After working one year as Lecturer, she joined the Division of Applied Mathematics and Process Control at IFP Energies nouvelles as a Research Assistant. Since 2005, she is an assistant professor in automatic control at the Department of Electrical and Chemical Engineering, UCB Lyon 1. Her research interests include state estimation problems, nonlinear observers, identification, and Lyapunov stability of nonlinear systems. She is also interested in control theory applications, especially in Biological systems and automotive applications.
\end{IEEEbiography}
\begin{IEEEbiography}[{\includegraphics[width=1in,height=1.25in,clip,keepaspectratio]{BioPhoto/Photo_Serres.jpg}}]%
{Ulysse Serres}
 received the Ph.D. degree in mathematics from the University of Burgundy, Dijon, France, in 2006. He is now Assistant Professor since 2009 at the Department of Processes and Electrical Engineering at the University Claude Bernard, Lyon, France
and researcher at the LAGEP Laboratory.
\end{IEEEbiography}
\fi

\end{document}